\documentclass{cmip-l} 
\usepackage{amsmath}
\usepackage{amssymb} 
\usepackage{amscd}
\usepackage{amsthm}
\usepackage[utf8]{inputenc}

\newtheorem{theorem}{Theorem}

\newtheorem{proposition}[theorem]{Proposition}

\newtheorem{corollary}[theorem]{Corollary}

\def\dim{{\mbox{dim}}}
\def\tr{{\mbox{tr}}}

\def\ker{{\mbox{Ker}}}     

\def\det{{\mbox{det}}}
\def\GL{{\mbox{GL}}}
\def\Hom {{\mbox{Hom}}}
\def\Ext {{\mbox{Ext}}}

\def\cali{{\mathcal I}}
 
\def\cala{{\mathcal A}} 
 
\def\cale{{\mathcal E}} 
\def\calk{{\mathcal K}} 
\def\call{{\mathcal L}}
\def\calm{{\mathcal M}}

\def\calh{{\mathcal H}}
\def\cals{{\mathcal S}}
 
\def\cale{{\mathcal E}}
 
\def\calf{{\mathcal F}}

\def\calw{{\mathcal W}}

\def\fracB{{\mathfrak B}}

\def\fracS{{\mathfrak {S}}}

\def\bbbone{\mbox{\rm 1\hspace {-.6em} l}}

\def\ug{{\mathbf u}}

\def\Tor{{\mbox{Tor}}}
\def\Algg{{\mbox{\bf Alg}}}

\numberwithin{equation}{section}

\begin{document}
\enlargethispage{3cm}

 \thispagestyle{empty}

\begin{center}
{\bf NONCOMMUTATIVE COORDINATE ALGEBRAS}
\end{center}

\vspace{0.3cm}

\begin{center}

Michel DUBOIS-VIOLETTE
\footnote{Laboratoire de Physique Th\'eorique, UMR 8627\\
Universit\'e Paris XI,
B\^atiment 210\\ F-91 405 Orsay Cedex\\
Michel.Dubois-Violette$@$u-psud.fr}\\
\end{center}
\vspace{0,5cm}
 \vspace{1cm}
 \begin{center}
 {\sl Dédié à Alain Connes}
 \end{center}
 
 \vspace{1cm}
 \begin{center}
 {\sl Il est bon de lire entre les lignes, cela fatigue moins les yeux}.\\
  Sacha Guitry
 \end{center}

 \begin{abstract}
We discuss the noncommutative generalizations of polynomial algebras which after appropriate completions can be used as coordinate algebras in various noncommutative settings, (noncommutative differential geometry, noncommutative algebraic geometry, etc.). These algebras have finite presentations and are completely characterized and classified by their (noncommutative) volume forms.
\end{abstract}

\tableofcontents

\section*{Introduction}
The universal skeletons for coordinate algebras in classical geometry (differential geometry, algebraic geometry, etc.) are polynomial algebras. The appropriate function algebras are obtained by completions with respect to the adapted topologies
 and either by gluing process or by taking quotient algebras.\\
Our aim here is to discuss the noncommutative generalizations of polynomial algebras which can be used similarily in various noncommutative setting, noncommutative differential geometry \cite{ac:1980}, \cite{ac:1986a}, \cite{ac:1994}, noncommutative algebraic geometry, etc. as well as in the applications in physics.\\

At the very beginning, one has to face the question of what class of algebras should we consider as generalization of the algebras of polynomial functions on finite-dimensional vector spaces. It seems clear that one must stay within the class of the $\mathbb N$-graded algebras which are connected, generated in degree 1 with finite presentations. There is a minimal choice which is the class of quadratic algebras which are Koszul (see below) of finite global dimension, which have polynomial growth and satisfy a version of Poincaré duality refered to as the Gorenstein property in \cite{art-sch:1987}. A bigger class is the class of regular algebras in the sense of \cite{art-sch:1987} which shall be refered to as the class of AS-regular algebras in the following. We shall consider here a bigger class in that we shall drop the condition of polynomial growth included in the AS-regularity condition. We shall refer to this bigger class of algebras as regular algebras. Although polynomial growth is a very natural condition for noncommutative coordinate algebras (and from the point of view of deformation theory), it turns out that for our analysis we do not need it and that by imposing polynomial growth one eliminates algebras which in spite of the fact that they have not an interpretation of (noncommutative) coordinate algebra are very interesting and are furthermore relevant for physics. Of course, at any stage one can restrict attention to the subclass of algebras with polynomial growth, (or which are quadratic, etc.). For global dimensions $D=2$ and $D=3$, the regular algebras are $N$-homogeneous and Koszul. We shall recall what this Koszul property means. This is a very desirable property that one can formulate only for $N$-homogeneous algebras for the moment (i.e. algebras with relations of degree $N$). This is why, for global dimensions $D\geq 4$ we shall impose $N$-homogeneity and Koszulity.\\

In the following we shall review various concepts and results. We shall in particular give a survey of the results of \cite{mdv:2005}, \cite{mdv:2007} in which we shall insist on the conceptual points and drop technical proofs.We shall illustrate the mains points by a lot of examples. A central result is that the algebras under consideration are completely specified by multilinear forms on finite-dimensional vector spaces. Given such an algebra, the corresponding multilinear form, which is unique up to a nonvanishing scale factor, plays the role of the (noncommutative) volume form. Furthermore isomorphic algebras correspond to multilinear forms which are in the same orbit of the corresponding linear group ($GL(g)$ for $g$ generators). The determination of the moduli space of these algebras is of course of mathematical interest by itself. Concerning physics, the classification of these algebras can become of great importance since in a noncommutative geometrical approach to the quantum theory of space and gravitation one should expect the occurrence at some approximation of a superposition of noncommutative geometries.\\

It is worth noticing that the results of \cite{mdv:2007} have been recently generalized to the quiver case in \cite{boc-sch-wem:2008}. The correspondence  between \cite{mdv:2007} and \cite{boc-sch-wem:2008} should read : multilinear forms or volumes $\leftrightarrow$ superpotentials.\\

Finally one should point out that this article is not only a survey but that it also contains new results and concepts.\\

Let us give some indications on the notations. Throughout the paper $\mathbb K$ denotes a field, all vector spaces and algebras are over $\mathbb K$, the dual of a vector space $E$ is denoted by $E^\ast$ and the symbol $\otimes$ denotes the tensor product over $\mathbb K$. Without other specifications, an algebra will always be an associative unital algebra. A graded algebra will be a $\mathbb N$-graded algebra $\cala=\oplus_{n\in \mathbb N} \cala_n$. Such a graded algebra is said to be connected whenever $\cala_0=\mathbb K\bbbone$. Given a $(r,s)$-matrix $A$, we denote by $A^t$ its transposed $(s,r)$-matrix. We use the Einstein summation convention of repeated up down indices in the formulas.

\section{Regular algebras}

The aim of this section is to make explicit the general class of algebras that we wish to investigate and to set up some notations.

\subsection{Graded algebras}

The algebras that we shall consider will be connected $\mathbb N$-graded algebras which are finitely generated in degree 1 and finitely presented with homogeneous relations of degrees $\geq 2$. These algebras are the objects of the category $\mathbf{GrAlg}$, the morphisms of this category being the homogeneous algebra homomorphisms of degree 0.\\
An algebra ${\cala}\in \mathbf{GrAlg}$ is of the form $\cala=A(E,R)=T(E)/[R]$ where $E=\cala_1$ is finite-dimensional, $R=\oplus_{n\geq 2} R_n$
 a finite-dimensional graded subspace of $T(E)$ such that (independence)
 \[
 R_n\cap [\oplus_{m<n} R_m]=\{0\}
 \]
 for any $n$ $(R_n=\{0\}$ for $n<2)$ and where $[F]$ denotes for any subset $F\subset T(E)$ the two-sided ideal generated by $F$. The graded vector space $R$ is the space of independent relations of $\cala$. \\
 By chosing a basis $(x^\lambda)_{\lambda\in \{1,\dots,g\}}$ of $E$ and an homogeneous basis $(f_\alpha)_{\alpha\in \{1,\dots,r\}}$ of $R$ one can also write
 \[
 \cala =\mathbb K\langle x^1,\dots,x^g\rangle/[f_1,\dots,f_r]
 \]
 where $f_\alpha \in E^{\otimes^{N_\alpha}}$, $N_\alpha\geq 2$. Notice that $r(=\dim R)$ is well defined (i.e. only depends on $\cala$).\\
 If $R$ is concentrated in degree $N\ (\geq 2)$ i.e. if $R\subset E^{\otimes^N}$ then $\cala$ will be said to be a $N$-homogeneous algebra. The $N$-homogeneous algebras form a full subcategory $\mathbf{H_N\Algg}$ of $\mathbf{GrAlg}$, \cite{ber:2001a}, \cite{ber-mdv-wam:2003}.
 
\subsection{Dimension}\label{Dim}

Let $\cala\in \mathbf{GrAlg}$ be as above so
$\cala=\mathbb K\langle x^1,\dots,x^g\rangle/[f_1,\dots,f_r]$ and one can define $M_{\alpha\lambda}\in E^{\otimes^{N_{\alpha}-1}}$ by setting
$f_\alpha=M_{\alpha\lambda} \otimes x^\lambda \in E^{\otimes^{N_\alpha}}$. Then the presentation of $\cala$ by generators and relations is equivalent to the exactness of the sequence of left $\cala$-modules \cite{art-sch:1987}
\begin{equation}
\cala^r \stackrel{M}{\rightarrow} \cala^g \stackrel{x}{\rightarrow} \cala \stackrel{\varepsilon}{\rightarrow} \mathbb K \rightarrow 0
\label{pr}
\end{equation}
where $M$ means right multiplication (in $\cala$) by the matrix $(M_{\alpha\lambda})$, $x$ means right multiplication by the column $(x^\lambda)$ and $\varepsilon$ is the projection onto $\cala_0=\mathbb K$. In more intrinsic notations the exact sequence (\ref{pr}) reads for $\cala=A(E,R)$ 
\begin{equation}
\cala\otimes R \rightarrow \cala\otimes E \stackrel{m}{\rightarrow} \cala\stackrel{\varepsilon}{\rightarrow} \mathbb K\rightarrow 0
\label{PR}
\end{equation}
where $m$ is the multiplication of $\cala$ and the first arrow is as in (\ref{pr}). The exact sequence (\ref{PR}) corresponding to the presentation of $\cala$ extends as a minimal projective resolution 
\[
\dots \rightarrow \cale_n \rightarrow \cale_{n-1}\rightarrow \cale_0\rightarrow \mathbb K \rightarrow 0
\]
of the left $\cala$-module $\mathbb K$ which is in fact a free resolution \cite{car:1958}
\begin{equation}
\dots \rightarrow \cala\otimes E_n\rightarrow \cala\otimes E_{n-1}\rightarrow \dots \rightarrow \cala\rightarrow \mathbb K \rightarrow 0
\label{mr}
\end{equation}
and it follows from the very definition of $\Ext_\cala(\mathbb K,\mathbb K)$ that one can make the identifications
\begin{equation}
E^\ast_n=\Ext^n_\cala(\mathbb K,\mathbb K)
\label{ext}
\end{equation}
which read $R^\ast=\Ext^2_\cala(\mathbb K,\mathbb K)$ and $E^\ast=\Ext^1_\cala(\mathbb K, \mathbb K)$ for $n=2$ and $n=1$. The Yoneda algebra $\Ext_\cala(\mathbb K,\mathbb K)$ is the cohomology of a graded differential algebra from which it follows that it carries a canonical $A_\infty$-structure \cite{kad:1980}, \cite{lu-pal-wu-zha:2006}. It turns out that one can reconstruct the graded algebra $\cala$ from the $A_\infty$-algebra $\Ext_\cala(\mathbb K, \mathbb K)$ \cite{kel:2001}, \cite{lu-pal-wu-zha:2006}. Thus the $A_\infty$-algebra $\Ext_\cala(\mathbb K, \mathbb K)$ is a natural dual of the graded algebra $\cala$. In the case of a $N$-homogeneous algebra $\cala$, there is another natural dual of $\cala$ which is its Koszul dual $\cala^!$ \cite{ber-mdv-wam:2003} (see below). In the case of a Koszul algebra these two notions are strongly connected and coincide in the quadratic case $(N=2)$, \cite{ber-mar:2006}.\\

The length of the resolution (\ref{mr}) is the projective dimension of the left module $\mathbb K$. It is classical \cite{car:1958},  \cite{art-tat-vdb:1990} that the left global dimension of $\cala$ (for $\cala\in \mathbf{GrAlg}$) coincides with the projective dimension of $\mathbb K$ as left module and that it also coincides with the right global dimension (and with the projective dimension of $\mathbb K$ as right module). Furthermore it has been shown recently \cite{ber:2005} that this dimension also coincides with the Hochschild dimension of $\cala$ in homology as well as in cohomology. So for an algebra $\cala\in \mathbf{GrAlg}$ there is a unique definition of the dimension from a homological point of view which will be refered to as {\sl its global dimension} in the sequel. In the following, we shall only consider algebras in $\mathbf{GrAlg}$ with finite global dimensions.\\
It is worth noticing here that there is another dimension for $\cala\in \mathbf{GrAlg}$ which is the Gelfand-Kirillov dimension but since in the following polynomial growth plays no role (and therefore will not be assumed) we shall only consider the global dimension for our general analysis.

\subsection{Poincaré duality}

We now assume that $\cala=A(E,R)\in \mathbf{GrAlg}$ is of finite global dimension $D$. The (free) resolution (\ref{mr}) of $\mathbb K$ reads then
\[
0\rightarrow \cala\otimes E_D\rightarrow \dots \rightarrow \cala\otimes E_1\rightarrow \cala \rightarrow \mathbb K \rightarrow 0
\]
with $E_1=\cala_1=E$.\\
By applying the functor $\Hom_\cala(\bullet,\cala)$ to the chain complex
\[
0\rightarrow \cala\otimes E_D\rightarrow \dots \rightarrow \cala\otimes E \rightarrow \cala\rightarrow 0
\]
of (free) left $\cala$-modules, one obtains a cochain complex $\cale'$
\[
0\rightarrow \cale'_0\rightarrow \cale'_1\rightarrow \dots \rightarrow \cale'_D\rightarrow 0
\]
of right $\cala$-modules. The cohomology $H(\cale')$ of this complex is by definition $\Ext_\cala(\mathbb K,\cala)$ that is one has
\begin{equation}
H^n(\cale')=\Ext^n_\cala(\mathbb K, \cala)
\label{extA}
\end{equation}
for any $n\in \mathbb N$.\\
By definition $\cala$ is said to be Gorenstein if one has
$\Ext^D_\cala(\mathbb K,\cala)=\mathbb K$ and $\Ext^n_\cala(\mathbb K, \cala)=0$
for $n\not= D$. This means that
\[
0\rightarrow \cale'_0\rightarrow \dots \rightarrow \cale'_D\rightarrow \mathbb K \rightarrow 0
\]
is a free resolution of $\mathbb K$ as right $\cala$-module. This resolution is then a minimal projective resolution of the right $\cala$-module $\mathbb K$ which implies the isomorphisms
\[
E^\ast_n \simeq E_{D-n}
\]
of vector spaces and therefore
\begin{equation}
\dim (E_n)=\dim (E_{D-n})
\label{DP}
\end{equation}
for $0\leq n\leq D$.\\
Thus the Gorenstein property is a variant of the Poincaré duality property.
 			
\subsection{Regularity}

Let $\cala=A(E,R)$ be a graded algebra of $\mathbf{GrAlg}$, $\cala$ will be said to be {\sl regular} if it is of finite global dimension, $g\ell\dim (\cala)=D<\infty$, and is Gorenstein. This definition of regularity is directly inspired from the one of \cite{art-sch:1987} which will be refered to as AS-regularity, the only difference is that we have dropped the condition of polynomial growth since we do not need it for the analysis in the sequel and since it would eliminate very interesting examples. \\
This is the class of algebras that we would like to analyse and we shall do it for low global dimensions $D=2$ and $D=3$. For higher global dimension, we shall restrict a little the class of algebras that we will consider. In order to understand this let us recall the following result \cite{ber-mar:2006}.

\begin{proposition}\label{R2-3}
Let $\cala$ be a regular algebra of global dimension $D$.\\
$\mathrm{(i)}$ If $D=2$ then $\cala$ is quadratic and Koszul.\\
$\mathrm{(ii)}$ If $D=3$ then $\cala$ is $N$-homogeneous with $N\geq 2$ and Koszul.
\end{proposition}
Thus for $D<4$, regularity implies $N$-homogeneity (with $N=2$ for $D=2$) and Koszulity. We shall explain later what is the Koszul property. This is a very desired property that one can formulate for the moment only for homogeneous algebras. This is why we shall restrict attention in the following to regular algebras which are $N$-homogeneous (with $N\geq 2$) and Koszul. In view of the above proposition this is not a restriction for regular algebras of global dimension $D=2$ and $D=3$ however one knows examples of regular algebras in global dimension 4 and more which are not homogeneous.

\section{Global dimension $D=2$}\label{globdim2}

This section is devoted to the description of the regular algebras of global dimension $D=2$.

\subsection{General results}

Let us use the notations of the beginning of $\S \ref{Dim}$ so let $\cala=\mathbb K\langle x^1,\dots,x^g \rangle/[f_1,\dots,f_r]$ and consider the exact sequence (\ref{pr}) corresponding to the presentation of $\cala$. The algebra $\cala$ has global dimension $D=2$ if and only if (\ref{pr}) extends as an exact sequence
\[
0\rightarrow \cala^r\stackrel{M}{\rightarrow} \cala^g \stackrel{x}{\rightarrow} \cala \stackrel{\varepsilon}{\rightarrow} \mathbb K \rightarrow 0
\]
i.e. as a free resolution of $\mathbb K$ of length $D=2$.\\
Assume now that $D=2$ and that $\cala$ is Gorenstein. Then the Gorenstein property implies that $r=1$, that degree $(M)$ = degree $(x)$ = 1 so $M=(B_{\rho\lambda}x^\rho)$ and that the matrix $(B_{\lambda\mu})\in M_g(\mathbb K)$ is invertible. The above free resolution of $\mathbb K$ reads then
\begin{equation}
0\rightarrow \cala \stackrel{x^tB}{\rightarrow} \cala^g \stackrel{x}{\rightarrow} \cala\stackrel{\varepsilon}{\rightarrow} \mathbb K \rightarrow 0
\label{res2}
\end{equation} 
with obvious notations.\\
Conversely, let $b$ be a nondegenerate bilinear form on $\mathbb K^g$ with matrix elements $B_{\lambda\mu}$ in the canonical basis and let $\cala$ be the (quadratic) algebra generated by $g$ generators $x^\lambda$ with relation $B_{\lambda\mu}x^\lambda x^\mu=0$, then $\cala$ is Gorenstein of global dimension $D=2$. One has the following theorem \cite{mdv:2007}.

\begin{theorem}\label{REG2}
Let $b$ be a nondegenerate bilinear form on $\mathbb K^g$ $(g\geq 2)$ with components $B_{\mu\nu}=b(e_\mu,e_\nu)$ in the canonical basis $(e_\lambda)$ of $\mathbb K^g$. Then the quadratic algebra $\cala$ generated by the elements $x^\lambda$ $(\lambda\in \{1,\dots,g\})$ with the relation $B_{\mu\nu}x^\mu x^\nu=0$ is regular of global dimension $D=2$. Conversely any regular algebra of global dimension $D=2$ is of the above kind for some $g\geq 2$ and some nondegenerate bilinear form $b$ on $\mathbb K^g$. Furthermore two such algebras $\cala$ and $\cala'$ are isomorphic if and only if $g=g'$ and $b'=b\circ L$ for some $L\in GL(g,\mathbb K)$.
\end{theorem}

The last part of this theorem is almost obvious and gives a description of the moduli space of the regular algebras of global dimension $D=2$.\\
The right action $b\mapsto b\circ L$ of the linear group on bilinear forms is a particular case of the right action of the linear group $GL(V)$ on multilinear forms on a vector space $V$ defined for a $n$-linear form $w$ by
\begin{equation}
w\circ L (v_1,\dots,v_n)=w(Lv_1,\dots,Lv_n)
\label{GLA}
\end{equation}
for any $v_k\in V$, $k\in \{1,\dots,n\}$.\\
For reasons which will become clear, the algebra $\cala$ (regular of global dimension $D=2$) associated to the nondegenerate bilinear form $b$ on $\mathbb K^g$ as in Theorem \ref{REG2} will be denoted $\cala(b,2)$ in the following.

\subsection{Poincaré series and polynomial growth}\label{Poi}

Let $\cala$ be a regular algebra of global dimension $D=2$. Then the exact sequence (\ref{res2}) splits as
\[
0\rightarrow \cala_{n-2}\stackrel{x^t B}{\rightarrow}\cala^g_{n-1}\stackrel{x}{\rightarrow} \cala_n\rightarrow 0
\]
for $n\not= 0$ with of course $\cala_0=\mathbb K$ and $\cala_n=0$ for $n<0$. It follows that the Poincaré series $P_\cala(t)$ of $\cala$ is given by
\begin{equation}
P_\cala (t)=\frac{1}{1-gt+t^2}
\label{poinca2}
\end{equation}
in view of the Euler-Poincaré formula.\\

For $g=2$ one has
\[
P_\cala(t)=\left(\frac{1}{1-t}\right)^2
\]
so $\cala$ has then polynomial growth (with $GK\dim=2$) while for $g>2$ one has
\[
P_\cala(t)=\frac{1}{(1-k^{-1}t)(1-kt)}
\]
with 
\[
k=\frac{1}{2}(g+\sqrt{g^2-4})>1
\]
so $\cala$ has then exponential growth.\\
Let us now discuss the case of the regular algebras of global dimension 2 with $g=2$ generators i.e. which have polynomial growth. In view of Theorem \ref{REG2} these algebras are classified by the $GL(2,\mathbb K)$-orbits of nondegenerate bilinear forms on $\mathbb K^2$. Assuming that $\mathbb K$ is algebraically closed, it is easy to classify these $GL(2,\mathbb K)$-orbits of nondegenerate bilinear forms according to the rank $\mathbf{rk}$ of their symmetric parts \cite{mdv-lau:1990} :\\

\noindent (0) $\mathbf{rk}=0$ - there is only one orbit which is the orbit of the bilinear form $b=\varepsilon$ with matrix of components
\[
B=\left(
\begin{array}{cc}
0 & -1\\
1 & 0
\end{array}
\right)
\]
which corresponds to the relations $x^1x^2-x^2x^1=0$ so $\cala$ is isomorphic to the polynomial algebra $\mathbb K[x^1,x^2]$, \\

\noindent (1) $\mathbf{rk}=1$ - there is only one orbit which is the orbit of the bilinear form $b$ with matrix of components
\[
B=\left(
\begin{array}{cc}
0 & -1\\
1 & 1
\end{array}
\right)
\]
which corresponds to the relations $x^1x^2-x^2x^1-(x^2)^2=0$,\\

\noindent (2) $\mathbf{rk}=2$ - there is a 1-parameter family of orbits which are the orbits of the bilinear forms $b=\varepsilon_q$ with matrices of components
\[
B=\left(
\begin{array}{cc}
0 & -1\\
q & 0
\end{array}
\right)
\]
for $q\in \mathbb K$ with $q^2-q\not= 0$ modulo $q\sim q^{-1}$ which corresponds to the relations $x^1x^2-qx^2x^1=0$.\\
The case (0) corresponds to the ordinary plane, the case (1) corresponds to the Jordanian plane and the cases (2) correspond to the Manin planes. One thus recovers the usual description of the algebras which are regular in the sense of \cite{art-sch:1987} i.e. AS-regular of global dimension 2, \cite{irv:1979}, 
\cite{art-sch:1987}.

\subsection{Hecke symmetries}\label{Hec}

Any linear mapping
\[
R:\mathbb K^g \otimes \mathbb K^g \rightarrow \mathbb K^g \otimes \mathbb K^g
\]
is characterized by its components $R^{\mu\nu}_{\lambda\rho}$ defined by
\[
R(e_\lambda \otimes e_\rho)=R^{\mu\nu}_{\lambda\rho} e_\mu \otimes e_\nu
\]
in the canonical basis $(e_\lambda)$ of $\mathbb K^g$.\\
Let $b$ be a nondegenerate bilinear form on $\mathbb K^g$ with components $B_{\lambda\rho}=b(e_\lambda,e_\rho)$ and let $K^{\mu\nu}$ be the components of a bilinear form on the dual vector space of $\mathbb K^g$ in the dual basis of $(e_\lambda)$. Define then the endomorphism $R$ of $\mathbb K^g \otimes \mathbb K^g$ by setting 
\begin{equation}
R^{\mu\nu}_{\lambda\rho}=\delta^\mu_\lambda \delta^\nu_\rho + K^{\mu\nu} B_{\lambda\rho}
\label{defR}
\end{equation}
for $\mu,\nu, \lambda,\rho \in \{1,\dots,g\}$. Assume now that the above $R$ defined by (\ref{defR}) satisfies the Yang-Baxter equation
\begin{equation}
(I\otimes R) (R\otimes I) (I\otimes R)=(R\otimes I)(I\otimes R)(R\otimes I)
\label{YB}
\end{equation}
on $(\mathbb K^g)^{\otimes^3}$ where $I$ denotes the identity mapping of $\mathbb K^g$ onto itself. One verifies that (\ref{YB}) is equivalent to

\begin{equation}
\left\{
\begin{array}{l}
KBK^tB^t + (1+\tr (KB^t)) \bbbone = 0 \\
K^t B^t KB + (1+ \tr (KB^t)) \bbbone = 0
\label{eqYB}
\end{array}
\right.
\end{equation}
where $K$ and $B$ are the matrices $(K^{\mu\nu})$ and $(B_{\lambda\rho})$ of $M_g(\mathbb K)$ and where the product is the matrix product. Equations (\ref{eqYB}) imply then that one has
\begin{equation}
(R-\bbbone)(R-(1+\tr (K B^t))\bbbone)=0
\label{Hk}
\end{equation}
which means that $R$ is a Hecke symmetry in the terminology of \cite{gur:1990}.\\

Given the nondegenerate bilinear for $b$, one can always solve (\ref{eqYB}).
For instance
\begin{equation}
K=qB^{-1}
\label{sol0}
\end{equation}
with $q\in \mathbb K$ such that
\begin{equation}
q +q^{-1}+\tr (B^{-1}B^t)=0
\label{q}
\end{equation}
is a solution of (\ref{eqYB}). The corresponding Hecke symmetries will be called {\sl the standard Hecke symmetries associated with} (the nondegenerate bilinear form) $b$ while more generally the Hecke symmetries associated with the solutions of (\ref{eqYB}) will be said to be {\sl associated with} $b$. There are generically two standard Hecke symmetries corresponding to the two roots of Equation (\ref{q}).\\
Notice that (\ref{eqYB}) implies that $K\not=0$ so if $R$ is a Hecke symmetry associated to $b$, the defining relation of $\cala(b,2)$ namely $B_{\mu\nu}x^\mu x^\nu=0$ is equivalent to the quadratic relations
\begin{equation}
x^\mu x^\nu=R^{\mu\nu}_{\lambda\rho} x^\lambda x^\rho
\label{Hcom}
\end{equation}
for $\mu,\nu\in \{1,\dots,g\}$.\\
In the case $g=2$ with $b=\varepsilon_q$ i.e.

\[
B=\left(\begin{array}{c c}
0 & -1\\
q & 0
\end{array}
\right)
\]
with $q\not=0$ which includes cases (0) and (2) of \S  
\ref{Poi}, one can take
\[
K=\left( \begin{array}{cc}
0 & 1\\
-p & 0
\end{array}
\right)
,\ p\in \mathbb K
\]
as solution of (\ref{eqYB}). Equation (\ref{Hk}) reads then
\[
(R-\bbbone)(R+pq)=0
\]
and for $p=q,R$ is a standard Hecke symmetry for $b=\varepsilon_q$. In the classical situation $q=1$, i.e. for $\cala=\mathbb K [x^1,x^2]$, both standard Hecke symmetries coincide and reduce to the flip
\[
x\otimes y \mapsto y\otimes x
\]
of $\mathbb K^2\otimes \mathbb K^2$.

\subsection{Actions of quantum groups}

There are quantum groups acting on the noncommutative planes corresponding to the regular algebras of global dimension $D=2$. For the Manin planes corresponding to the $\cala(\varepsilon_q,2)$ these are the quantum groups $SL_q(2)$, $GL_q(2)$ and $GL_{p,q}(2)$ \cite{man:1987}, \cite{man:1988}.\\
For the noncommutative plane corresponding to $\cala(b,2)$ where $b$ is a nondegenerate bilinear form on $\mathbb K^g$, the generalization of $SL_q(2)$ is the quantum group of the nondegenerate bilinear form $b$, \cite{mdv-lau:1990}.
Let us recall the definition of this object. Let $\calh(b)$ be the unital associative algebra generated by the $g^2$ elements $u^\mu_\nu$ ($\mu,\nu \in \{1,\dots, g\}$) with the relations
\begin{equation}
B_{\lambda\rho} u^\lambda_\mu u^\rho_\nu=B_{\mu\nu} \bbbone
\label{Binv1}
\end{equation}
and
\begin{equation}
B^{\mu\nu} u^\lambda_\mu u^\rho_\nu=B^{\lambda\rho} \bbbone
\label{Binv2}
\end{equation}
where the $B^{\mu\nu}$ are the matrix elements of the inverse matrix $B^{-1}$ of the matrix $B$ of the components $B_{\mu\nu}=b(e_\mu,e_\nu)$ of $b$. One verifies easily that there is a unique structure of Hopf algebra on $\calh(b)$ with coproduct $\Delta$, counit $\varepsilon$ and antipode $S$ such that

\begin{eqnarray}
\Delta (u^\mu_\nu)& = &u^\mu_\rho \otimes u^\rho_\nu\label{cp}\\
\varepsilon(u^\mu_\nu) & = & \delta^\mu_\nu \label{cu}\\
S(u^\mu_\nu) & = & B^{\mu\lambda}B_{\rho\nu} u^\rho_\lambda\label{anti}
\end{eqnarray}
the product and the unit being the original ones on $\calh(b)$.\\
There is a canonical algebra-homomorphism $\Delta_L:\cala(b,2)\rightarrow \calh(b)\otimes \cala(b,2)$ such that
\[
\Delta_L(x^\lambda)=u^\lambda_\mu \otimes x^\mu
\]
for $\lambda\in \{1,\dots,g\}$. This equips $\cala(b,2)$ with a structure of $\calh(b)$-comodule. The dual object of $\calh(b)$ is the quantum group of the nondegenerate bilinear form $b$. The analysis of the category of representations of this quantum group has been done in \cite{bic:2003b}.
 To the coaction $\Delta_L$ of $\calh(b)$ on $\cala(b,2)$ corresponds an action of this quantum group on the noncommutative plane corresponding to $\cala(b,2)$.\\
The (quadratic) homogeneous part of the relations (\ref{Binv1}) and (\ref{Binv2}) reads
\begin{equation}
u^\mu_\alpha u^\nu_\beta R^{\alpha\beta}_{\lambda\rho}=R^{\mu\nu}_{\alpha\beta} u^\alpha_\lambda u^\beta_\rho
\label{Bhom}
\end{equation}
where $R$ is a standard Hecke symmetry of $b$. In fact (\ref{Bhom}) together with (\ref{cp}) and (\ref{cu}) define a bialgebra with counit for any $R$. In the case where $R$ is a standard Hecke symmetry then $B^{\mu\nu}B_{\rho\lambda} u^\lambda_\mu u^\rho_\nu$ is in the center and the Hopf algebra $\calh(b)$ corresponding to the quantum group of $b$ is the quotient of the bialgebra by the ideal generated by the element
\[
B^{\mu\nu} B_{\rho\lambda}u^\lambda_\mu u^\rho_\nu-g \bbbone
\]
of the center. In fact $\calh(b)$ is a quotient of a bigger Hopf algebra associated with the homogeneous relations (\ref{Bhom}) which is the generalization of the Hopf algebra corresponding to $GL_q(2)$ in the case $b=\varepsilon_q$, ($g=2$). More generally if $R$ is arbitrary an Hecke symmetry associated with $b$, there is a Hopf algebra associated with the quadratic relations (\ref{Bhom}) which coacts on $\cala(b,2)$ and corresponds to the generalization of $GL_{p,q}(2)$.

\section{Global dimension $D=3$}\label{globdim3}

In this section we shall analyse regular algebras of global dimension $D=3$ and describe some representative examples. For global dimensions $D\geq 3$ what replace the bilinear forms of the global dimension $D=2$ (last section) are multilinear forms so we start this section with a discussion on multilinear forms.

\subsection{Multilinear forms}

Let $V$ be a vector space with $\dim (V)\geq 2$, $Q$ be an element of the linear group $GL(V)$ and $m$ be an integer with $m\geq 2$. Then a $m$-linear form $w$ on $V$ (i.e. a linear form on $V^{\otimes^m}$) will be said to be $Q$-{\sl cyclic} if one has 
\begin{equation}
w(X_1,\dots,X_m)=w(QX_m,X_1,\dots,X_{m-1})
\label{Qcycl}
\end{equation}
for any $X_1,\dots,X_m\in V$.\\
Let $w$ be $Q$-cyclic then one has 
\[
w(X_1,\dots,X_m)=w(QX_k,\dots,QX_m,X_1,\dots, X_{k-1})
\]
for $1\leq k\leq m$ so in particular one has
\[
w(X_1,\dots,X_m)=w(QX_1,\dots,QX_m)
\]
for any $X_1,\dots,X_m\in V$ which also reads $w=w\circ Q$ and means that $w$ is invariant by $Q$.\\
Let now $w$ be an artibrary $Q$-invariant $m$-linear form on $V$, then the $m$-linear form $\pi_Q(w)$ on $V$ defined by
\[
\pi_Q(w) (X_1,\dots,X_m)=\frac{1}{m}\sum^m_{k=1} w(QX_k,\dots,QX_m,X_1,\dots,X_{k-1})
\]
for any $X_1,\dots,X_m\in V$ is $Q$-cyclic and this defines a projection $\pi_Q$ of the space of $Q$-invariant $m$-linear forms onto the space of $Q$-cyclic $m$-linear forms on $V$. This projection is $GL(V)$-equivariant in the sense that if $w$ is $Q$-invariant (resp. $Q$-cyclic) then $w\circ L$ is $L^{-1}QL$-invariant (resp. $L^{-1}QL$-cyclic) for any $L\in\GL(V)$.\\

The $m$-linear form $w$ on $V$ will be said to be {\sl preregular} if it satisfies the following conditions (i) and (ii) :\\
(i) $w(X,X_1,\dots,X_{m-1})=0$ for any $X_1,\dots,X_{m-1}\in V$ implies $X=0$,\\
(ii) there is a $Q_w\in GL(V)$ such that $w$ is $Q_w$-cyclic.\\
Condition (i) implies that $Q_w$ is unique under (ii) and Condition (ii) and (i) imply that $w$ satisfies the following condition (i') which is stronger than (i) : \\
(i') $w$ $(X_1,\dots, X_k,X,X_{k+1},\dots, X_{m-1})=0$ for any $X_1,\dots, X_{m-1} \in V$ implies 
\[
X=0,\  \mbox{for any}\  k\in \{0,\dots, m-1\}.
\]
A $m$-linear form $w$ on $V$ satisfying (i') will be said to be 1-{\sl site-nondegenerate}.\\
The set of preregular $m$-linear forms on $V$ is invariant by the action of $GL(V)$ and one has
\begin{equation}
Q_{w\circ L}=L^{-1}Q_wL
\label{Prinv}
\end{equation}
for any preregular $m$-linear form $w$ on $V$.\\
A bilinear form $b$ on $\mathbb K^g$ is preregular if and only if it is nondegenerate ; one has then $Q_b=(B^{-1})^tB$ where $B$ is the matrix of components of $b$.\\

The condition of preregularity will be involved throughout the paper. We now introduce a stronger condition which is involved specifically in the description of the regular algebras of global dimension $D=3$. Let $N$ be an integer with $N\geq 2$, then a $(N+1)$-linear form $w$ on $V$ will be said to be 3-{\sl regular} if it is preregular and satisfies the following condition (iii) :\\
(iii)\hspace{1cm}  If $L_0$ and $L_1$ are endomorphisms of $V$ satisfying
\[
w(L_0X_0,X_1,X_2,\dots,X_N)=w(X_0,L_1 X_1,X_2,\dots,X_N)
\]
for any $X_0,\dots,X_n\in V$, then $L_0=L_1=k\bbbone$ for some $k\in \mathbb K$.\\
The set of 3-regular $(N+1)$-linear forms is also invariant by $GL(V)$.\\
Condition (iii) is a sort of 2-sites nondegeneracy condition. Consider the stronger condition (iii') :\\
(iii')\hspace{1cm} $\sum_i w (Y_i,Z_i,X_1,\dots,X_{N-1})=0$ for any $X_1,\dots,X_{N-1}\in V$ implies 
\[
\sum_i Y_i\otimes Z_i=0.
\]
 It is clear that (iii') $\Rightarrow$ (iii), however it is a strictly stronger condition. For instance let $\varepsilon$ be the completely antisymmetric $(N+1)$-linear form on $\mathbb K^{N+1}$ with $\varepsilon(e_0,\dots,e_N)=1$. Then $\varepsilon$ is 3-regular but one has
\[
\varepsilon (Y,Z,X_1,\dots,X_{N-1})+\varepsilon (Z,Y,X_1\dots, X_{N-1})=0
\]
 identically and this does not imply $Y\otimes Z + Z\otimes Y=0$.

\subsection{General results for $D=3$}

Let $w$ be a preregular $(N+1)$-linear form on $\mathbb K^g$ with components $W_{{\lambda_0}\dots {\lambda_N}}=w(e_{\lambda_0},\dots,e_{\lambda_N})$ in the canonical basis ($e_\lambda$) of $\mathbb K^g$ and let $\cala(w,N)$ be the $N$-homogeneous algebra generated by the $g$ elements $x^\lambda$ ($\lambda\in\{1,\dots,g\}$) with the $g$ relations
\begin{equation}
W_{\lambda\lambda_1\dots\lambda_N}x^{\lambda_1}\dots x^{\lambda_N}=0
\label{Re3}
\end{equation}
for $\lambda\in \{1,\dots,g\}$. In other words one has
$\cala(w,N)=A(E,R)$ with $E=\oplus_\lambda \mathbb K x^\lambda$
and $R=\sum_\lambda \mathbb K W_{\lambda\lambda_1\dots \lambda_N} x^{\lambda_1}\otimes \dots \otimes x^{\lambda_N}$. Condition (i) implies that $\dim(R)=g$ that is that the latter sum is direct and that the relations (\ref{Re3}) are independent. \\

Let us now use again the notations of the beginning of \S \ref{Dim} so let $\cala\in \mathbf{GrAlg}$ with $\cala=\mathbb K\langle x^1,\dots,x^y\rangle/[f_1,\dots,f_r]$ and consider the exact sequence (\ref{pr}) corresponding to the presentation of $\cala$. Then $\cala$ has global dimension $D=3$ if and only if (\ref{pr}) extends as an exact sequence
\[
0\rightarrow \cala^s\rightarrow \cala^r \stackrel{M}{\rightarrow} \cala^g \stackrel{x}{\rightarrow}\cala \stackrel{\varepsilon}{\rightarrow} \mathbb K \rightarrow 0
\]
i.e. as a free resolution of $\mathbb K$ of length $D=3$. Assume now that $\cala$ is regular. Then the Gorenstein property (Poincaré duality) implies immediately that $r=g$, that $s=1$, that the above resolution reads with an appropriate choice of the relations $f_\lambda$
\begin{equation}
0\rightarrow \cala\stackrel{x^t}{\rightarrow} \cala^g \stackrel{M}{\rightarrow} \cala^g\stackrel{x}{\rightarrow} \cala \stackrel{\varepsilon}{\rightarrow} \mathbb K \rightarrow 0
\label{res3}
\end{equation}
and that $w=x^\lambda\otimes f_\lambda$ is homogeneous, say of degree $N+1$, and is preregular \cite{art-sch:1987}. So $\cala=\cala(w,N)$ as above. In fact one has the following theorem 
\cite{mdv:2007}.

\begin{theorem}\label{REG3}
Let $\cala$ be a regular algebra of global dimension $D=3$. Then $\cala=\cala(w,N)$ for some $N\geq 2$, some $g\geq 2$ and some 3-regular $(N+1)$-linear form $w$ on $\mathbb K^g$.
\end{theorem}
The Poincaré series of $\cala=\cala(w,N)$ as in the above theorem (i.e. $\cala$ regular with $D=3$) is given by
\begin{equation}
P_\cala(t)=\frac{1}{1-gt+gt^N-t^{N+1}}
\label{poinca3}
\end{equation}
in view of (\ref{res3}).\\

If one compares this theorem with Theorem \ref{REG2} for $D=2$, one sees that there are two missing items : first there is no converse of the statement in Theorem \ref{REG3} and second there is no characterization of the isomorphism classes. Concerning the first point, it was conjectured in \cite{mdv:2007} that given a 3-regular $(N+1)$-linear form $w$ on $\mathbb K^g$ then $\cala(w,N)$ is a regular algebra with $D=3$, but unfortunately this is wrong and we shall give counter-examples (see below). This means that one has to find some slightly stronger condition than 3-regularity for $w$ (for $D=3$).\\
Concerning the second point the following result holds (independently of the regularity of the algebras) \cite{mdv:2007}.

\begin{proposition}\label{Iso3}
Let $w$ be a 3-regular $(N+1)$-linear form on $\mathbb K^g$ and let $w'$ be a 3-regular $(N'+1)$-linear form on $\mathbb K^{g'}$. Then $\cala(w,N)$ and $\cala(w',N')$ are isomorphic if and only if $g'=g$, $N'=N$ and $w'=w\circ L$ for some $L\in GL(g,\mathbb K)$.
\end{proposition}
The conditions $g'=g$ and $N'=N$ are clear but the 3-regularity is really involved in the proof of this proposition (see in \cite{mdv:2007}).\\

Following \cite{art-sch:1987} one deduces from (\ref{poinca3}) that a regular algebra of global dimension $D=3$ has polynomial growth if and only if $g=3$ and $N=2$ or $g=2$ and $N=3$; Otherwise it has exponential growth (for $g\geq 2$ and $N\geq 2$).

\subsection{Examples and counter-examples}

All AS-regular algebras of global dimension $D=3$ give of course examples and our notations $w, M, Q_w$ come from \cite{art-sch:1987}. In fact, the classification of the regular algebras of global dimension $D=3$ with polynomial growth is based on the possible Jordan decompositions of the corresponding $Q_w$'s. Let us give some representative examples.\\
(a) {\sl The 3-dimensional Sklyanin algebra} \cite{ode-fei:1989}, \cite{ode:2002}. This is the algebra $\cala$ generated by 3 elements $x, y, z$ with relations
\begin{equation}
\left\{
\begin{array}{l}
xy-qyx=pz^2\\
yz-qzy=px^2\\
zx-qxz=py^2
\end{array}
\right.
\label{Sk3}
\end{equation}
where $p,q\in \mathbb K$ with $(p,q)\not= (0,0)$ and $(p^3+1,q^3+1)\not= (0,0)$.\\
This algebra is AS-regular with $D=3$. One has $\cala=\cala(w,2)$ with 
\begin{equation}
\begin{array}{lll}
w & = & x \otimes y \otimes z + y \otimes z \otimes x + z \otimes x \otimes y\\
& - & q(x \otimes z\otimes y + y\otimes x\otimes z + z \otimes y \otimes x)\\
& - & p (x\otimes x \otimes x + y \otimes y \otimes y + z \otimes z \otimes z)
\end{array}
\label{wSk3}
\end{equation}
where we have identified the 3-linear form $w$ on $\mathbb K^3$ with the corresponding element of $(\mathbb K^{3\ast})^{\otimes^3}$. One verifies that $w$ is 3-regular and one has
\begin{equation}
Q_w=\bbbone
\label{QSk3}
\end{equation} 
for the corresponding element of $GL(3,\mathbb K)$.\\

\noindent (b) {\sl The $q$-deformed 3-dimensional polynomial algebra}.
This is the algebra $\cala$ generated by 3 elements $x,y,z$ with relations

\begin{equation}
\left\{
\begin{array}{l}
xy=qcyx\\
yz=qazy\\
zx=qbxz
\end{array}
\right.
\label{qdef3}
\end{equation}
with $q,a,b,c\in \mathbb K$, $abc=1$ and $q\not=0$. This algebra is AS-regular with $D=3$ and one has $\cala=\cala(w,2)$ with 
\begin{equation}
w=bx\otimes y\otimes z+cy\otimes z\otimes x+az\otimes x\otimes y-q(abx\otimes z\otimes y+bcy\otimes x\otimes z+caz\otimes y\otimes x)
\label{wqdef3}
\end{equation}
with the same conventions as above. One verifies that $w$ is 3-regular and one has
\begin{equation}
Q_w=\left(
\begin{array}{ccc}
b/c & 0 & 0\\
0 & c/a & 0\\
0 & 0 & a/b
\end{array}
\right )
\label{Qqdef3}
\end{equation}

\noindent (c) {\sl Type E quadratic AS-algebra \cite{art-sch:1987}}. 
This is the algebra $\cala$ generated by 3 elements $x,y,z$ with relations
\begin{equation}
\left\{
\begin{array}{l}
x^2+\zeta^{-1}yz + \zeta\  zy =0\\
y^2+\zeta^{-4}zx + \zeta^4 xz=0\\
z^2+\zeta^{-7} xy + \zeta^7yx=0
\end{array}
\right.
\label{E}
\end{equation}
where $\zeta\in \mathbb K$ is a primitive 9th root of 1, $\zeta^9=1$.\\
This algebra is AS-regular with $D=3$ and $\cala=\cala(w,2)$ with 
\begin{equation}
\begin{array}{lll}
w & = & x\otimes z\otimes x + y\otimes x\otimes y + z\otimes y\otimes z\\
& + & \zeta\  z\otimes x\otimes x +\zeta^{-1}x\otimes x\otimes z\\
& + & \zeta^4 x\otimes y\otimes y + \zeta^{-4}y\otimes y \otimes x\\
& + & \zeta^7 y\otimes z\otimes z + \zeta^{-7} z\otimes z \otimes y
\end{array}
\label{wE}
\end{equation}
which defines a 3-regular 3-linear form on $\mathbb K^3$.\\
One has
\begin{equation}
Q_w=\left(
\begin{array}{ccc}
\zeta & 0 & 0\\
0 & \zeta^4 & 0\\
0 & 0 & \zeta^7
\end{array}
\right)
\label{QE}
\end{equation}
for the corresponding element of $GL(3,\mathbb K)$.\\
It is worth noticing here that the algebras of Case (a) and Case (b) are deformations of the polynomial algebra $\mathbb K[x,y,z]$ while this is not the case here. In fact the algebra with relations (\ref{E}) is quite rigid.\\

\noindent (d) {\sl Counter-example to the converse of Theorem \ref{REG3}}.
Let $\cala$ be the algebra generated by 3 elements $x,y,z$ with relations 
\begin{equation}
\left\{
\begin{array}{l}
x^2+yz=0\\
y^2+zx=0\\
xy=0
\label{Cex}
\end{array}
\right.
\end{equation}
Then $\cala=\cala(w,2)$ where the 3-linear form $w$ on $\mathbb K^3$ is given by
\begin{equation}
w=x\otimes x\otimes x + y \otimes y \otimes y + x\otimes y\otimes z + y\otimes z \otimes x + z\otimes x \otimes y
\label{wCex}
\end{equation}
with the same conventions as before. One verifies that $w$ is again 3-regular and one has $Q_w=\bbbone$. However $\cala$ is not regular of global dimension $D=3$. Indeed the candidate for (\ref{res3}) is
\[
0\rightarrow \cala \stackrel{x^t}{\rightarrow} \cala^3\stackrel{M}{\rightarrow} \cala^3\stackrel{x}{\rightarrow} \cala \stackrel{\varepsilon}{\rightarrow} \mathbb K \rightarrow 0
\]
with $x^t=(x,y,z)$ and
\[
M=\left(
\begin{array}{ccc}
x & 0 & y\\
z & y & 0\\
0 & x & 0\\
\end{array}
\right)
\]
but this complex is {\sl not} exact in second position : One has
\[
(yz,0,0)\in \ker (M)
\]
while $(yz,0,0)$ is not in the image of $x^t$.\\
This algebra is discussed in \cite{art-sch:1987} and there is a similar one which is cubic with 2 generators.\\

\noindent (e) {\sl The Yang-Mills algebra \cite{ac-mdv:2002b}}.
The Yang-Mills algebra is the cubic algebra $\cala$ generated by $g$ elements $\nabla_\lambda$ ($\lambda\in \{1,\dots,g\}$) with relations
\begin{equation}
g^{\lambda\mu}[\nabla_\lambda,[\nabla_\mu,\nabla_\nu]]=0
\label{YM}
\end{equation}
for $\nu\in \{1,\dots,g\}$, where the $g^{\lambda\mu}$ are the components of a symmetric nondegenerate bilinear form on $\mathbb K^g$. The use here of covariant instead of contravariant notations-conventions has a physical origin. This algebra is regular of global dimension $D=3$.
One has $\cala=\cala(w,3)$ where $w$ is the 4-linear form on $\mathbb K^g$ with components
\begin{equation}
W^{\alpha_1\alpha_2\alpha_3\alpha_4}=g^{\alpha_1\alpha_2}g^{\alpha_3\alpha_4}+g^{\alpha_2\alpha_3}g^{\alpha_4\alpha_1}-2g^{\alpha_1\alpha_3}g^{\alpha_2\alpha_4}
\label{wYM}
\end{equation}
for $\alpha_k\in \{1,\dots,g\}$. This 4-linear form on $\mathbb K^g$ is 3-regular, in fact it satisfies the strong condition (iii'), and one has $Q_w=\bbbone$.\\

\noindent (f) {\sl The super Yang-Mills algebra \cite{ac-mdv:2007}}.
There is a ``super" version of the Yang-Mills algebra which is the cubic algebra $\tilde \cala$ generated by $g$ elements $S_\lambda\ (\lambda\in \{1,\dots,g\})$ with relations
\begin{equation}
g^{\lambda\mu}[S_\lambda,[S_\mu,S_\nu]_+]=0
\label{SYM}
\end{equation}
for $\nu\in\{1,\dots,g\}$, where the $g^{\lambda\mu}$ are as above and $[A,B]_+=AB+BA$.\\
This algebra is again regular of global dimension 3 and $\tilde \cala=\cala(\tilde w,3)$ where $\tilde w$ is the 4-linear form on $\mathbb K^g$ with components
\begin{equation}
\tilde W^{\alpha_1\alpha_2\alpha_3\alpha_4}=g^{\alpha_2\alpha_3}g^{\alpha_4\alpha_1}-g^{\alpha_1\alpha_2}g^{\alpha_3\alpha_4}
\label{wSYM}
\end{equation}
for $\alpha_k\in \{1,\dots,g\}$. This $\tilde w$ is 3-regular (and satisfies (iii')) and $Q_{\tilde w}=-\bbbone$. Notice that the equations (\ref{SYM}) are equivalent to 
\begin{equation}
[S_\lambda,g^{\mu\nu}S_\mu S_\nu]=0
\label{CS2}
\end{equation}
i.e. to the fact that $g^{\mu\nu}S_\mu S_\nu$ is central.\\

Before leaving this section, it is worth noticing that the Yang-Mills algebra is by its very definition the universal enveloping algebra of a graded Lie algebra. In the case $g=2$ this is an AS-regular algebra considered in \cite{art-sch:1987} which is the universal enveloping algebra of the graded 3-dimensional Lie algebra with basis ($\nabla_1,\nabla_2$) in degree 1 and $C$ in degree 2 with Lie bracket defined by
\[
[\nabla_1,\nabla_2]=C,[\nabla_1,C]=0, [\nabla_2,C]=0
\]
In the case $g>2$, the Yang-Mills algebra has exponential growth.\\
Similar considerations apply to the super Yang-Mills where the above Lie algebra is replaced by a super Lie algebra.

\section{Homogeneous algebras}

The aim of this section is to describe properties of $N$-homogeneous algebras and to introduce and discuss the Koszul property \cite{ber:2001a}, \cite{ber-mdv-wam:2003}.

\subsection{Koszul duality}

Let $\cala\in {\mathbf H_N}\Algg$ be a $N$-homogeneous algebra, that is $\cala=A(E,R)$ with $R\subset E^{\otimes^N}$. One defines the (Koszul) {\sl dual} $\cala^!$ of $\cala$ to be the $N$-homogeneous algebra 
\begin{equation}
\cala^!=\cala(E^\ast, R^\perp)
\label{Kd}
\end{equation} 
where $R^\perp \subset E^{\ast\otimes^N}=(E^{\otimes^N})^\ast$ is the annihilation of $R$, i.e. the subspace 
\[
R^\perp = \{\omega\in (E^{\otimes^N})^\ast\vert \omega(x)=0,\  \forall x\in R\}
\]
 of $(E^{\otimes^N})^\ast$ identified with $E^{\ast\otimes^N}$
(there a canonical identification of $(E^{\otimes^N})^\ast$ with $E^{\ast\otimes^N}$ since $E$ is finite-dimensional). One has canonically
\begin{equation}
(\cala^!)^!=\cala
\label{Kdd}
\end{equation}
and to any morphism $f:\cala\rightarrow \cala'=A(E',R')$ of $\mathbf{H_N}\Algg$ corresponds a morphism $f^!:\cala'{^!}\rightarrow \cala^!$ which is induced by
the transposed of the restriction $f \restriction E: E\rightarrow E'$ of $f$ to $E$. The correspondence $(\cala\mapsto \cala^!,f\mapsto f^!)$ defines a contravariant involutive functor $((f^!)^!=f)$.

\subsection{The Koszul $N$-complex $K(\cala)$}.

Let $\cala=A(E,R)$ be a $N$-homogeneous algebra with dual $\cala^!=\oplus_n \cala^!_n$ and consider the dual vector spaces $\cala^{!\ast}_n$ of the $\cala^!_n$. One has

\begin{equation}
\left\{
\begin{array}{l}
\cala^{!\ast}_n=E^{\otimes^n}\ \ \text{for}\ \ n<N\\
\cala^{!\ast}_n=\cap_{r+s=n-N} E^{\otimes^r}\otimes R \otimes E^{\otimes^s}\ \ \text{for}\ \ n\geq N
\end{array}
\right.
\label{dstar}
\end{equation}
so that for any $n\in \mathbb N$ one has $\cala^{!\ast}_n\subset E^{\otimes^n}$. Let us define then the sequence of homomorphisms of (free) left $\cala$-modules
\begin{equation}
\dots \stackrel{d}{\rightarrow} \cala \otimes \cala^{!\ast}_{n+1} \stackrel{d}{\rightarrow} \cala\otimes \cala^{!\ast}_n\stackrel{d}{\rightarrow}\dots \stackrel{d}{\rightarrow} \cala\rightarrow 0
\label{KNC}
\end{equation}
where $d:\cala\otimes \cala^{!\ast}_{n+1}\rightarrow \cala\otimes \cala^{!\ast}_n$ is induced by the map
\[
a\otimes (e_0\otimes e_1 \otimes\dots \otimes e_n)\mapsto ae_0 \otimes (e_1\otimes \dots \otimes e_n)
\]
of $\cala\otimes E^{\otimes^{n-1}}$ into $\cala\otimes E^{\otimes^n}$. Then one has
\begin{equation}
d^N=0
\label{NCd}
\end{equation}
since $\cala^{!\ast}_n\subset R\otimes E^{\otimes^{n-N}}$ for $n\geq N$. Thus (\ref{KNC}) defines a $N$-complex which will be refered to as {\sl the Koszul $N$-complex of $\cala$} and denoted by $K(\cala)$.\\
As for any $N$complex \cite{mdv:1998a} one obtains from $K(\cala)$ of family $C_{p,r}(K(\cala))$ of ordinary complexes, called {\sl the contractions of} $K(\cala)$, by putting together alternatively $p$ and $N-p$ arrows $d$ of $K(\cala)$. The complex $C_{p,r}(K(\cala))$ is defined as
\begin{equation}
\dots \stackrel{d^{N-p}}{\rightarrow} \cala\otimes \cala^{!\ast}_{Nk+r} \stackrel{d^p}{\rightarrow} \cala\otimes \cala^{!\ast}_{Nk-p+r} \stackrel{d^{N-p}}{\rightarrow} \cala\otimes \cala^{!\ast}_{N(k-1)+r}\stackrel{d^p}{\rightarrow} \dots
\label{Cpr}
\end{equation}
for $0\leq r<p\leq N-1$ \cite{ber-mdv-wam:2003}, (one verifies that all such complexes are exhausted by these couples $(p,r)$). For the homology of these complexes one has the following result \cite{ber-mdv-wam:2003}.

\begin{proposition}\label{CPR}
Let $\cala=A(E,R)$ be a $N$-homogeneous algebra with $N\geq 3$. Assume that $(p,r)$ is distinct from $(N-1,0)$ and that $C_{p,r}(K(\cala))$ is exact at degree $k=1$. Then $R=0$ or $R=E^{\otimes^N}$.
\end{proposition}

In other words except for $C_{N-1,0}(K(\cala))$ a nontrivial acyclicity for $C_{p,r}(K(\cala))$ leads to the trivial algebras $\cala=T(E)$ or $\cala=T(E)^!$.

\subsection{Koszul complexes and Koszul property}

The last proposition points out the complex $C_{N-1,0}(K(\cala))$ which will be denoted by $\calk(\cala,\mathbb K)$ and refered to as {\sl the Koszul complex of $\cala$}. It coincides with the Koszul complex originality introduced in \cite{ber:2001a} without mention to the $N$-complex $K(\cala)$. Of course for a quadratic algebra $\cala$, i.e. for $N=2$, one has $K(\cala)=\calk (\cala,\mathbb K)$ and this coincides with the definition of \cite{pri:1970} (see also \cite{man:1987}, \cite{man:1988}).\\

A $N$-homogeneous algebra $\cala$ will be said to be a {\sl Koszul algebra} whenever its Koszul complex $\calk(\cala,\mathbb K)$ is acyclic in positive degrees, (i.e. $H_n(\calk (\cala,\mathbb K))=0$ for $n\geq 1$). This is the generalization given in 
\cite{ber:2001a} of the definition of \cite{pri:1970} to $N$-homogeneous algebras. There are very good reasons explained in \cite{ber:2001a} for this generalization. We content ourselves here to observe that, among the contractions of $K(\cala)$, the Koszul complex $\calk(\cala,\mathbb K)$ is distinguished by the fact that it terminates as a projective resolution of $\mathbb K$. Indeed the presentation of $\cala=A(E,R)$ is equivalent to the exactness of the sequence
\[
\cala\otimes R \stackrel{d^{N-1}}{\rightarrow} \cala\otimes E \stackrel{d}{\rightarrow} \cala \stackrel{\varepsilon}{\rightarrow} \mathbb K \rightarrow 0
\]
as observed before and, on the other hand one has $\cala^{!\ast}_1=E$ and $\cala^{!\ast}_N=R$ so $\calk(\cala,\mathbb K)$ terminates as
\[
\dots \stackrel{d}{\rightarrow} \cala\otimes R \stackrel{d^{N-1}}{\rightarrow} \cala\otimes E \stackrel{d}{\rightarrow} \cala\rightarrow 0
\]
thus if $\cala$ is a Koszul algebra, one has a free resolution of $\mathbb K$ which is then in fact a minimal projective resolution of the trivial left $\cala$-module $\mathbb K$ given by
\begin{equation}
\calk(\cala,\mathbb K)\stackrel{\varepsilon}{\rightarrow} \mathbb K \rightarrow 0
\label{KRes}
\end{equation}
which is refered to as the Koszul resolution of (the left $\cala$-module) $\mathbb K$.\\
Notice that if $\cala$ is a regular algebra of global dimension 2 (resp. 3) then (\ref{res2}) (resp. (\ref{res3})) are the Koszul resolutions of $\mathbb K$ (with a slight abuse of language) so that $\cala$ is then a Koszul algebra as announced in Proposition \ref{R2-3}. One has the following result \cite{mdv-pop:2002}.

\begin{proposition}\label{PSKN}
Let $\cala$ be a Koszul $N$-homogeneous algebra. One has
\[
P_\cala(t)Q_\cala(t)=1
\]
where the series $Q_\cala(t)$ is defined by
\[
Q_\cala(t)=\sum_{n\in \mathbb N} (\dim (\cala^!_{Nn})t^{Nn} - \dim (\cala^!_{Nn+1})t^{Nn+1})
\]
and where $P_\cala(t)=\sum_n \dim(\cala_n)t^n$
is the Poincaré series of $\cala$.
\end{proposition}
In fact the Koszul $N$-complex splits into sub-$N$complexes for the total degree
\[
K(\cala)=\oplus K^{(n)}(\cala)
\]
which induces a splitting of the Koszul complex into finite-dimensional subcomplexes
\[
\calk(\cala,\mathbb K)=\oplus \calk^{(n)}(\cala,\mathbb K)
\]
and the proposition follows from the Euler-Poincaré formula applied to each components.\\
Notice that in the quadratic case, one has $Q_\cala(t)=P_{\cala^!}(-t)$. \\

If $\cala$ is a Koszul $N$-homogeneous algebra, one has clearly
\begin{equation}
\cala^!_{Nn}\simeq \Ext^{2n}_\cala(\mathbb K,\mathbb K),\cala^!_{Nn+1}\simeq \Ext^{2n+1}_\cala(\mathbb K,\mathbb K)
\label{ExtN}
\end{equation}
and therefore by setting
\begin{equation}
Y_\cala(t)=\sum_{n\in \mathbb N} (\dim (\Ext^{2n}_\cala(\mathbb K,\mathbb K))t^{Nn}-\dim(\Ext^{2n+1}_\cala(\mathbb K,\mathbb K))t^{Nn+1})
\label{PYN}
\end{equation}
one has $P_\cala(t)Y_\cala(t)=1$. In \cite{kri:2007} it is shown that, conversely if $\cala$ is a $N$-homogeneous algebra such that one has
\begin{equation}
P_\cala(t)Y_\cala(t)=1
\label{NumK}
\end{equation}
then $\cala$ is Koszul. This gives an interesting numerical criterion for Koszulity which has to be compared with the fact that there are $N$-homogeneous algebras $\cala$ satisfying $P_\cala(t)Q_\cala(t)=1$ which are not Koszul, (of course then (\ref{ExtN}) do not hold).\\

In (\ref{KNC}) the factors $\cala$ are considered as left $\cala$-modules. By considering $\cala$ as a right $\cala$-module and by exchanging the factors, one obtains a $N$-complex $\tilde K(\cala)$ of right $\cala$-modules.
\begin{equation}
\dots\stackrel{\tilde d}{\rightarrow} \cala^{!\ast}_{n+1} \otimes \cala \stackrel{\tilde d}{\rightarrow} \cala^{!\ast}_n\otimes \cala \stackrel{\tilde d}{\rightarrow}\dots \stackrel{\tilde d}{\rightarrow} \cala\rightarrow 0
\label{RNC}
\end{equation}
where $\tilde d:\cala^{!\ast}_{n+1}\otimes \cala\rightarrow \cala^{!\ast}_n\otimes \cala$ is induced by the mapping $(e_1\otimes \dots\otimes e_{n+1})\otimes a \mapsto (e_1\otimes \dots \otimes e_n)\otimes e_{n+1}a$ of $E^{\otimes^{n+1}}\otimes \cala$ into $E^{\otimes^n}\otimes \cala$. The fact that $\tilde d^N=0$ follows from $\cala^{!\ast}_N\subset E^{\otimes^{n-N}}\otimes R$ for $n\geq N$. Let us consider the sequences $(L,R)$
\begin{equation}
\dots \stackrel{d_L,d_R}{\rightarrow} \cala\otimes \cala^{!\ast}_{n+1} \otimes \cala \stackrel{d_L,d_R}{\rightarrow}\cala \otimes \cala^{!\ast}_n \otimes \cala\rightarrow \dots \stackrel{d_L,d_R}{\rightarrow} \cala\otimes\cala\rightarrow 0
\label{BBNC}
\end{equation}
where $d_L=d\otimes I$ and $d_R=I\otimes \tilde d$, $I$ being the identity mapping of $\cala$ onto itself. One has $d^N_L=d^N_R=0$ and $d_L$ and $d_R$ are homomorphisms of $(\cala,\cala)$-bimodule, i.e. of left $\cala\otimes \cala^{opp}$-modules. The two $N$-differentials $d_L$ and $d_R$ commute so one has
\[
(d_L-d_R)\left(\sum^{N-1}_{p=0}d^p_L d^{N-p-1}_R\right)=\left(\sum^{N-1}_{p=0}d^p_L d^{N-p-1}_R\right)(d_L-d_R)=d^N_L-d^N_R=0.
\]
It follows that one defines a complex of free $\cala\otimes \cala^{opp}$-modules $\calk(\cala,\cala)$ by setting 
\begin{equation}
\left\{
\begin{array}{l}
\calk_{2m}(\cala,\cala)=\cala\otimes \cala^{!\ast}_{Nm}\otimes \cala\\
\\
\calk_{2m+1}(\cala,\cala)=\cala\otimes \cala^{!\ast}_{Nm+1}\otimes \cala
\end{array}
\right.
\label{KCB1}
\end{equation}
with differential $\delta'$ defined by
\begin{equation}
\left\{
\begin{array}{l}
\delta'=d_L-d_R:\calk_{2m+1}(\cala,\cala)\rightarrow \calk_{2m}(\cala,\cala)\\
\\
\delta'=\sum^{N-1}_{p=0}d^p_Ld^{N-p-1}_R:\calk_{2(m+1)}(\cala,\cala)\rightarrow \calk_{2m+1}(\cala,\cala)
\end{array}
\right.
\label{KCB2}
\end{equation}
which will be refered to as {\sl the bimodule Koszul complex of} $\cala$.\\
It turns out that $\calk(\cala,\cala)$ is acyclic in positive degrees if and only if $\calk(\cala,\mathbb K)$ is acyclic in positive degrees that is if and only if $\cala$ is a Koszul algebra. On the other hand one has the obvious exact sequence of bimodules
\[
\cala\otimes E\otimes \cala \stackrel{\delta'}{\rightarrow} \cala\otimes \cala \stackrel{m}{\rightarrow} \cala \rightarrow 0
\]
where $m$ denotes the product of $\cala$. This means that $H_0(\calk(\cala,\cala))=\cala$ and therefore whenever $\cala$ is Koszul one has a free resolution
\[
\calk(\cala,\cala)\stackrel{m}{\rightarrow}\cala\rightarrow 0
\]
of the left $\cala\otimes \cala^{opp}$-module $\cala$ which is a minimal projective resolution of $\cala$ and will be refered to as the Koszul resolution of $\cala$.

\subsection{Small complex and Poincaré duality for Koszul algebras}

Let $\cala$ be a $N$-homogeneous Koszul algebra and let $\calm$ be a $(\cala,\cala)$-bimodule considered as a right $\cala\otimes\cala^{opp}$-module. Then by interpreting the Hochschild homology $H(\cala,\calm)$ of $\cala$ with values in $\calm$ as $\Tor^{\cala\otimes \cala^{opp}}(\calm,\cala)$ \cite{car-eil:1973}, one sees that the homology of the complex $\calm \otimes_{\cala\otimes \cala^{opp}}\calk(\cala,\cala)$ is the $\calm$-valued Hochschild homology of $\cala$. We shall refer to this latter complex as the {\sl small Hochschild complex} of the Koszul algebra $\cala$ with coefficients in $\calm$ and denote it by $\cals(\cala,\calm)$. It reads
\begin{equation}
\dots \stackrel{\delta}{\rightarrow} \calm \otimes \cala^{!\ast}_{N(m+1)}\stackrel{\delta}{\rightarrow}\calm \otimes \cala^{!\ast}_{Nm+1}\stackrel{\delta}{\rightarrow}\calm\otimes \cala^{!\ast}_{Nm}\stackrel{\delta}{\rightarrow} \dots
\label{SKC}
\end{equation}
where $\delta$ is obtained from $\delta'$ by applying the factors $d_L$ to the right of $\calm$ and the factors $d_R$ to the left of $\calm$.\\

By construction the lengths of the complexes $\calk(\cala,\mathbb K)$ and $\calk(\cala,\cala)$ coincide. Assume that $\cala$ is a Koszul algebra, then this implies that the projective dimension of the trivial $\cala$-module $\mathbb K$ coincides with the Hochschild dimension of $\cala$ which is a particular case of the general result of 
\cite{ber:2005}.\\

The Koszul complex $\calk(\cala,\mathbb K)$ is a chain complex since its differential is of degree -1, the same is true for $\calk(\cala,\cala)$. By applying the functor $\Hom_\cala(\bullet, \cala)$ to the chain complex of free left $\cala$-modules $\calk(\cala,\mathbb K)$ one obtains the cochain complex $\call(\cala,\mathbb K)$ of free right $\cala$-modules
\[
0\rightarrow \call^0(\cala,\mathbb K)\rightarrow \dots \rightarrow \call^n(\cala,\mathbb K)\rightarrow \dots
\]
where $\call^n(\cala,\mathbb K)= \Hom_\cala(\calk_n(\cala,\mathbb K),\cala)$. 
Assume that $\cala$ is Koszul of global dimension $D$. Then $\call^n(\cala,\mathbb K)=0$ for $n>D$ and $\cala$ is Gorenstein if and only if $H^n(\call(\cala,\mathbb K))=0$ for $n<D$ and $H^D(\call(\cala,\mathbb K))=\mathbb K$. When $\cala$ is Koszul of global dimension $D$ and Gorenstein, this implies a precise form of the Poincaré duality between the Hochschild homology and the Hochschild cohomology of $\cala$, \cite{ber-mar:2006}, 
\cite{vdb:1998}, \cite{vdb:2002}. In the case of a regular algebra $\cala=\cala(w,N)$ of global dimension 3, it reads for an $\cala$-bimodule $\calm$
\begin{equation}
H_n(\cala,\calm)=H^{3-n}(\cala,\calm)
\label{HPd3}
\end{equation}
for $0\leq n\leq 3$ when $Q_w=\bbbone$, (when $Q_w\not = \bbbone$ it induces an automorphism $\sigma_w$ of $\cala$ and one has to twist by $\sigma_w$ the left multiplication of $\calm$ by $\cala$ on the right-hand side of (\ref{HPd3})).\\

The complex $\call(\cala, \mathbb K)$ is also a contraction of a natural $N$-complex $L(\cala)$. This $N$-complex $L(\cala)$ is the cochain $N$-complex of free right $\cala$-modules  obtained by applying the functor $\Hom_\cala(\bullet,\cala)$ to the Koszul $N$-complex $K(\cala)$ (which is a chain $N$-complex of free left $\cala$-modules). The right $\cala$-module $L^n(\cala)$ identifies canonically with $\cala^!_n\otimes \cala$ while the $N$-differential of $L(\cala)$ is then the left multiplication by $x^\ast_\lambda\otimes x^\lambda$ in $\cala^! \otimes \cala$ where $(x^\ast_\lambda)$ is the dual basis of $(x^\lambda)$ ($E=\oplus_\lambda \mathbb K x^\lambda$). One has $\call(\cala,\mathbb K)=C_{1,0}(L(\cala))$, i.e. $\call^0(\cala,\mathbb K)=\cala=L^0(\cala),\call^1(\cala,\mathbb K)=L^N(\cala)$, etc.

\subsection{Examples of Koszul algebras}
All regular algebras of global dimensions $D=2$ and $D=3$ are Koszul so in particular the examples of regular algebras of Sections 2 and 3 are examples of Koszul algebras.
We shall describe regular Koszul algebras of higher global dimension $D$ in Section 5. Let us give here some examples of Koszul algebras which are not generically regular.\\

\noindent (a) {\sl Koszul duals of quadratic algebras}.
It is well known and not hard to show that if $\cala$ is a quadratic algebra, then its Koszul dual $\cala^!$ is Koszul if and only if $\cala$ is Koszul. Even if $\cala$ is regular, $\cala^!$ is generically not regular.\\
For instance the exterior algebra $\wedge \mathbb K^g$ is the Koszul dual of the algebra of polynomial functions on $\mathbb K^g$ which is regular and Koszul of global dimension $g$, however $\wedge\mathbb K^g$ is not of finite global dimension.\\
It is worth noticing here that if $\cala$ is a $N$-homogeneous algebra with $N>2$, then the Koszulity of $\cala$ does not imply the Koszulity of its Koszul dual $\cala^!$, (this is due to the jumps in degrees in the Koszul resolution). For instance the Koszul dual $\cala^!$ of the Yang-Mills algebra $\cala$ (\S 3.3, example (e)) is such that $P_{\cala^!}(t)Q_{\cala^!}(t)\not=1$ (by direct computation) so it is not Koszul in view of Proposition \ref{PSKN}.\\

\noindent (b) {\sl Degenerate bilinear form} \cite{ber:2008}. In the following $b$ is a bilinear form on $\mathbb K^g$ with $g\geq 2$, $B=(B_{\mu\nu})$ is the matrix of components $B_{\mu\nu}=b(e_\mu,e_\nu)$ of $b$ in the canonical basis of $\mathbb K^g$ and $\cala=\cala(b,2)$ is the quadratic algebra generated by $g$ elements $x^\lambda$ with the relation
\[
B_{\mu\nu} x^\mu x^\nu=0
\]
i.e. we generalize the notation of Section 2 to cases where $b$ can be degenerate. In 
\cite{ber:2008} one finds the following results (Propositions 5.4 and 5.5 in \cite{ber:2008}) which contains Theorem \ref{REG2}.

\begin{proposition}\label{Roland}
Assume that $b\not=0$, then $\cala=\cala(b,2)$ has the following properties :\\
1) $\cala$ is Koszul,\\
2) $\cala$ has global dimension $D=2$ except in the case where $b$ is symmetric of rank 1 in which case $D=\infty$,\\
3) $\cala$ is Gorenstein if and only if $b$ is nondegenerate.
\end{proposition}

Thus for $b$ degenerate one has a lot of examples of Koszul algebras which are not regular. In \cite{ber:2008} there is a similar statement for $N$-homogeneous algebras with one relation ($r=1$) which although slightly more involved permits the construction of examples (see e.g. Example (c) in the next section \S 5.3).\\

(c) {\sl The self-duality algebra} \cite{ac-mdv:2002b}. In the case $g=4$ and $g^{\lambda\mu}=\delta^{\lambda\mu}$, the Yang-Mills algebra (Example (e) in \S\ 3.3) admits the 2 nontrivial quotients $\cala^{(+)}$ and $\cala^{(-)}$ where $\cala^{(\varepsilon)}$ $(\varepsilon=\pm)$ is the quadratic algebra generated by the 4 elements $\nabla_\lambda$ ($\lambda\in\{1,2,3,4\}$) with relations
\begin{equation}
[\nabla_4,\nabla_k]=\varepsilon[\nabla_\ell, \nabla_m]
\label{SD}
\end{equation}
for any cyclic permutation $(k,\ell,m)$ of (1,2,3). Let us fix $\varepsilon=+$ and call $\cala^{(+)}$ the self-duality algebra (the study of $\cala^{(-)}$ is similar). In 
\cite{ac-mdv:2002b} it was shown that this algebra is Koszul of global dimension $D=2$ and that the Koszul resolution reads
\begin{equation}
0\rightarrow (\cala^{(+)})^3\rightarrow (\cala^{(+)})^4\rightarrow \cala^{(+)}\stackrel{\varepsilon}{\rightarrow}\mathbb K \rightarrow 0
\label{resSD}
\end{equation}
from which it follows that
\begin{equation}
P_{\cala^{(+)}}(t)=\frac{1}{(1-t)(1-3t)}
\label{poincaSD}
\end{equation}
so $\cala^{(+)}$ has exponential growth and is not Gorenstein.\\
It follows from the definition that $\cala^{(+)}$ is the universal enveloping algebra of a Lie algebra which is the semi-direct product of the free Lie algebra $L(\nabla_1,\nabla_2,\nabla_3)$ by the derivation $\delta$ given by
\begin{equation}
\delta(\nabla_k)=[\nabla_\ell, \nabla_m]
\label{DerL}
\end{equation}
for any cyclic permutation $(k,\ell,m)$ of (1,2,3). Formula (\ref{poincaSD}) as well as all the above properties of $\cala^{(+)}$ follow also directly from this structure.\\

\noindent (d) {\sl The super self-duality algebra} \cite{ac-mdv:2007}. In a similar way as in the last example, for $g=4$ and $g^{\lambda\mu}=\delta^{\lambda\mu}$, the super Yang-Mills algebra (Example (f) in \S\ 3.3) admits the 2 nontrivial quotients $\tilde \cala^{(+)}$ and $\tilde\cala^{(-)}$ where $\tilde\cala^{(\varepsilon)}$ ($\varepsilon=\pm$) is the quadratic algebra generated by the 4 elements $S_\lambda$ ($\lambda\in \{1,2,3,4\}$) with relations 
\begin{equation}
i[S_4,S_k]_+=\varepsilon[S_\ell,S_m]
\label{SSD}
\end{equation}
for any cyclic permutation $(k,\ell,m)$ of (1,2,3). Let us fix $\varepsilon=+$ and call $\tilde\cala^{(+)}$ the super self-duality algebra. This algebra is again a Koszul algebra of global dimension 2 which is not Gorenstein and has Poincaré series given by
\begin{equation}
P_{\tilde A^{(+)}}(t)=\frac{1}{(1-t)(1-3t)}
\label{PSSD}
\end{equation}
so has also exponential growth. This algebra has direct relations with the 4-dimensional Sklyanin algebra (see in \cite{ac-mdv:2007}).

\section{Arbitrary global dimension $D$}

In the previous sections, we have seen that the regular algebras of global dimensions $D=2$ and $D=3$ are $N$-homogeneous (with $N=2$ for $D=2$) and Koszul. This very desirable property permits to write explicit canonical resolutions. On the other hand one can formulate for the moment this Koszul property only for $N$-homogeneous algebras. This is why in this section we shall restrict attention to Koszul homogeneous algebras and our aim is then to formulate the generalization of Theorem \ref{REG3} for arbitrary global dimension $D$. Notice however that for global dimensions $D\geq 4$, regularity does not imply $N$-homogeneity. It is worth mentioning here that for $D=4$ the AS-regular algebras, i.e. the regular algebras with polynomial growth, have been recently classified \cite{lu-pal-wu-zha:2004}.\\

We shall need a class of $N$-homogeneous algebras associated with preregular multilinear forms that we now describe.

\subsection{Homogeneous algebras associated to multilinear forms}

In this subsection $m$ and $N$ are integers with $m\geq N\geq 2$ and $w$ is a preregular $m$-linear form on $\mathbb K^g$ ($g\geq 2)$ with components $W_{\lambda_1\dots\lambda_m}=w(e_{\lambda_1},\dots,e_{\lambda_m})$ in the canonical basis $(e_\lambda)$ of $\mathbb K^g$. Let $\cala=\cala(w,N)$ be the $N$-homogeneous algebra generated by the elements $x^\lambda$ ($\lambda\in \{1,\dots,g\}$) with relation
\begin{equation}
W_{\lambda_1\dots\lambda_{m-N}\mu_1\dots \mu_N}x^{\mu_1}\dots x^{\mu_N}=0
\label{ReD}
\end{equation}
for $\lambda_k\in \{1,\dots,g\}$. Thus one has $\cala=A(E,R)$ with $E=\oplus_\lambda\mathbb K x^\lambda$ and  
\[
R=\sum_{\lambda_k} \mathbb K W_{\lambda_1\dots \lambda_{m-N}}\mu_1\dots \mu_Nx^{\mu_1}\otimes \dots\otimes x^{\mu_N}\subset E^{\otimes^N}.
\]
Notice that this generalizes the definitions of Section \ref{globdim2} (which is the case $m=N=2$) and Section \ref{globdim3} (which is the case $m=N+1$).\\

Let us define the subspaces $\calw_n\subset E^{\otimes^n}$ for $m\geq n\geq 0$ by
\begin{equation}
\left\{
\begin{array}{l}
\calw_n=E^{\otimes^n}\ \ \ \text{for}\ \ \ N-1\geq n\geq 0\\
\calw_n=\sum_{\lambda_k}\mathbb K W_{\lambda_1\dots \lambda_{m-n}\mu_1\dots\mu_n}x^{\mu_1}\otimes \dots \otimes x^{\mu_n}\ \ \ \text{for}\ \ \ m\geq n\geq N
\end{array}
\right.
\label{SubK}
\end{equation}
so in particular $\calw_1=E$ and $\calw_N=R$. The twisted cyclicity of $w$ (property (ii) of \S 3.1) and (\ref{dstar}) implies the following proposition.

\begin{proposition}\label{SubNC}
The sequence 
\begin{equation}
0\rightarrow \cala\otimes \calw_m\stackrel{d}{\rightarrow} \cala\otimes \calw_{m-1}\stackrel{d}{\rightarrow} \dots \stackrel{d}{\rightarrow} \cala\rightarrow 0
\label{SubNK}
\end{equation}
is a sub-$N$-complex of the Koszul $N$-complex $K(\cala)$ of $\cala$.
\end{proposition}
In fact one has $\calw_n \subset \cala^{!\ast}_n$ and $d(\cala\otimes \calw_{n+1})\subset \cala\otimes \calw_n$. In particular one has $\calw_m=\mathbb K w\subset \cala^{!\ast}_m$ so $w$ is a linear form on $\cala^!_m$. We define then the linear form $\omega_w$ on the algebra $\cala^!$ by setting
\begin{equation}
\omega_w =w\circ p_m
\label{omegaw}
\end{equation}
where $p_m:\cala^!\rightarrow \cala^!_m$ is the canonical projection onto the degree $m$. With $E=\oplus_\lambda \mathbb K x^\lambda$, $w$ is canonically a $m$-linear form on $E^\ast$ and $Q_w$ an element of $GL(E^\ast)$. With these identifications one has the following theorem \cite{mdv:2007}.

\begin{theorem}\label{MOD}
The element $Q_w$ of $GL(E^\ast)$ induces an automorphism $\sigma_w$ of the $N$-homogeneous algebra $\cala^!=A(E^\ast, R^\perp)$ and one has
\begin{equation}
\omega_w(xy)=\omega_w(\sigma_w(y)x)
\label{pmod}
\end{equation}
for any $x,y \in \cala^!$. The subset of $\cala^!$
\[
\cali=\{ y\in \cala^!\vert \omega_w(xy)=0,\ \ \forall x\in \cala^!\}
\]
is a two-sided ideal of $\cala^!$ and the quotient algebra $\calf(w,N)=\cala^!/\cali$ equipped with the linear form induced by $\omega_w$ is a graded Frobenius algebra.
\end{theorem}
To prove this theorem, one first verifies by using the $Q_w$-invariance of $w$ that one has $Q^{\otimes^N}_w R^\perp \subset R^\perp$ which implies the existence of $\sigma_w$. Then (\ref{pmod}) is just a translation of the $Q_w$-cyclicity of $w$. By definition $\cali$ is a left ideal and (\ref{pmod}) implies that it is also a right ideal. The quotient $\calf=\cala^!/\cali$ is a finite-dimensional graded algebra and the pairing induced by $(x,y)\mapsto \omega_w(xy)$ is nondegenerate and is a Frobenius pairing on $\calf$.

\begin{corollary}\label{AUT}
Considered as an element of $GL(E)$, the transposed $Q^t_w=Q^w$ of $Q_w$ induces an automorphism $\sigma^w$ of the $N$-homogeneous algebra $\cala=A(E,R)$.
\end{corollary}

Let us end this subsection by noting that, at this level of generality and for $N=2$ (i.e. in the quadratic case), the multilinear form $w$ induces a (twisted) noncommutative $m$-form for $\cala$. For this let $^w\cala$ be the $(\cala,\cala)$-bimodule which coincides with $\cala$ as right $\cala$-module and is such that the structure of left $\cala$-module is given by the left multiplication by $(-1)^{(m-1)n}(\sigma^w)^{-1}(a)$ for $a\in \cala_n$. One has the following result \cite{mdv:2007}.

\begin{proposition}\label{FORM}
In the case $N=2$ that is for $\cala=\cala(w,2)$, $\bbbone \otimes w$ is canonically a non trivial $^w\cala$-valued Hochschild $m$-cycle on $\cala$.
\end{proposition}
 In this statement, $\bbbone$ is interpreted as an element of $^w\cala$ while $w\in E^{\otimes^m}$ is interpreted as an element of $\cala^{\otimes^m}$ ($E=\cala_1\subset \cala$) so that $\bbbone \otimes w$ is a $^w\cala$-valued Hochschild $m$-chain.
 
\subsection{General results for Koszul-Gorenstein algebras}

For the $N$-homo\-geneous algebras which are Koszul of finite global dimension $D$ and which are Gorenstein, (a particular class of regular algebras if $D\geq 4$), one has the following theorem \cite{mdv:2007}.

\begin{theorem}\label{KGD}
Let $\cala$ be a $N$-homogeneous algebra which is Koszul of finite global dimension $D$ and Gorenstein. Then $\cala=\cala(w,N)$ for some preregular $m$-linear form on $\mathbb K^g$ for some $g$. If $N\geq 3$ then $m=Np+1$ and $D=2p+1$ for some $p\geq 1$ while for $N=2$ one has $m=D$.
\end{theorem}

For the proof we refer to \cite{mdv:2007}.\\

Under the assumptions of Theorem \ref{KGD} the Koszul resolution of the trivial left $\cala$-module $\mathbb K$ reads
\[
0\rightarrow \cala\otimes \calw_m \stackrel{d}{\rightarrow} \cala \otimes \calw_{m-1} \stackrel{d^{N-1}}{\rightarrow} \dots \stackrel{d}{\rightarrow}\cala\otimes \calw_N \stackrel{d^{N-1}}{\rightarrow}\cala\otimes E\stackrel{d}{\rightarrow}\cala \stackrel{\varepsilon}{\rightarrow}\mathbb K \rightarrow 0
\]
or, by setting
\begin{equation}
\left\{
\begin{array}{l}
\nu_N(2k)=Nk\\
\nu_N(2k+1)=Nk+1
\end{array}
\right.
\label{nuN}
\end{equation}
for $k\in \mathbb N$,
\begin{equation}
0\rightarrow \cala\otimes \calw_{\nu_N(D)}\stackrel{d'}{\rightarrow} \dots \stackrel{d'}{\rightarrow} \cala\otimes \calw_{\nu_N(k)} \stackrel{d'}{\rightarrow} \cala\otimes \calw_{\nu_N(k-1)}\dots \stackrel{d'}{\rightarrow}\cala\stackrel{\varepsilon}{\rightarrow}\mathbb K\rightarrow 0
\label{KResD}
\end{equation}
where $d'$ is defined by
\begin{equation}
\left\{
\begin{array}{l}
d'=d^{N-1}:\cala\otimes \calw_{\nu_N(2k)}\rightarrow \cala\otimes \calw_{\nu_N(2k-1)}\\
d'=d : \cala\otimes \calw_{\nu_N(2k+1)}\rightarrow \cala\otimes \calw_{\nu_N(2k)}
\end{array}
\right.
\label{diffW}
\end{equation}
for $k\in \mathbb N$.\\
Notice that one has
\begin{equation}
\dim(\calw_{\nu_N(k)})=\dim (\calw_{\nu_N(D-k)})
\label{PoinDual}
\end{equation}
for $0\leq k\leq D$. In particular $\cala\otimes \calw_{\nu_N(D)}=\cala\otimes w$ so one sees that $\bbbone \otimes w$ is the generator of the top module of the Koszul resolution which again corresponds to the interpretation of $\bbbone\otimes w$ as a volume form.\\

It is worth noticing here that it has been already shown in \cite{bon-pol:1994} that the quadratic algebras which are Koszul and regular are determined by multilinear form ($D$-linear for global dimension $D$) which correspond to volume forms in this noncommutative setting.\\

Let us come back on a more general situation. Assume that $D$ and $N$ are given integers with $D\geq 2$ and $N\geq 2$ and that $N=2$ whenever $D$ is an even integer. Let then $w$ be a preregular $m$-linear form on $\mathbb K^g$ with $m=D$ for $N=2$ and $m=Np+1$ for $D=2p+1$ and consider the $N$-homogeneous algebra $\cala=\cala(w,N)$. The complex
\begin{equation}
0\rightarrow \cala \otimes \calw_{\nu_N(D)}\stackrel{d'}{\rightarrow} \dots \stackrel{d'}{\rightarrow} \cala\otimes \calw_{\nu_N(k)}\stackrel{d'}{\rightarrow} \dots \stackrel{d'}{\rightarrow} \cala\rightarrow 0
\label{CWD}
\end{equation}
is still well defined, with $\nu_N$ as in (\ref{nuN}) and $d'$ as in (\ref{diffW}), and is a subcomplex of the Koszul complex $\calk(\cala,\mathbb K)$ of $\cala$ in view of Proposition \ref{SubNC}. It is clear that if this complex is acyclic in positive degree, it coincides with the Koszul complex of $\cala$ and that $\cala$ is then Koszul of global dimension $D$ and Gorenstein. Thus as remarked in \cite{boc-sch-wem:2008} one has the following result which gives a sort of converse of Theorem \ref{KGD}.

\begin{proposition}\label{CKGD}
Let $\cala=\cala(w,N)$ be as above then $\cala$ is Koszul of global dimension $D$ and Gorenstein if and only if the complex $\mathrm{(\ref{CWD})}$ is acyclic in positive degrees.
\end{proposition}

A weaker assumption on the complex (\ref{CWD}) is to assume that it coincides with the Koszul complex. In the case where $D=3$, one has the following proposition \cite{mdv:2007}.

\begin{proposition}\label{Eq3R}
Let $w$ be a preregular $(N+1)$-linear form on $\mathbb K^g$ and let $\cala=\cala(w,N)$ then the following conditions are equivalent.\\

\noindent $\mathrm{(a)}$ $\cala^{!\ast}_{N+1}=\mathbb K w$.\\

\noindent $\mathrm{(b)}$ The complex $\mathrm{(\ref{CWD})}$ coincides with the Koszul complex $\calk(\cala,\mathbb K)$ of $\cala$.\\

\noindent  $\mathrm{(c)}$ $w$ is 3-regular.\\
\end{proposition}

Let us consider $\cala=\cala(w,N)$ with Koszul dual $\cala^!=\oplus_n \cala^!_n$ and let us define the graded algebra
\begin{equation}
\cala'=\cala'(w,N)=\oplus_n \cala'_n
\label{Aprime1}
\end{equation}
to be $\cala^!$ for $N=2$ and to be defined for $N>2$ by
\begin{equation}
\cala'_n=\cala^!_{\nu_N(n)}
\label{Aprime2}
\end{equation}
for $n\in \mathbb N$ with product $(x,y)\mapsto x\bullet y$ defined by
\begin{equation}
x\bullet y = \pi(xy)
\label{Aprime3}
\end{equation}
where $\pi:\cala^!\rightarrow \cala'$ is the canonical projection of $\cala^!$ onto $\cala'=\oplus_n \cala^!_{\nu_N(n)}\subset \cala^!$ defined by setting $\pi(\cala^!_k)=0$ whenever $k$ is not in $\nu_N(\mathbb N)$. Thus this product is defined for two homogeneous elements $x$ and $y$ by 
\[
xy=0
\]
whenever $x$ and $y$ are both of odd degree and 
\[
xy = \text{product in}\ \cala^!
\]
otherwise. It is clear that this product is associative. One has the following result.

\begin{theorem}\label{APR}
Assume that $D,N$ and $w$ are as above that is $N=2$ for $D$ even and $w$ is a preregular $m$-linear form on $\mathbb K^g$ with $m=D$ for $N=2$ and $m=Np+1$ for $D=2p+1$. Then the following conditions are equivalent.\\

\noindent  $\mathrm{(a)}$ $\cala'(w,N)$ equipped with the linear form induced by $\omega_w$ is a Frobenius algebra.\\

\noindent  $\mathrm{(b)}$ The complex $\mathrm{(\ref{CWD})}$ coincides with the Koszul complex $\calk(\cala,\mathbb K)$ of $\cala(w,N)$.\\
\end{theorem}

\noindent \underbar{Proof}. 
The proof of this proposition is almost tautological since conditions (a) and (b) are both equivalent to $\calw_{\nu_N(n)}=\cala^{!\ast}_{\nu_N(n)}={\cala'_n}^\ast$ for $n\in \mathbb N$. $\square$\\

This is of course directly inspired by \cite{ber-mar:2006} and implies Theorem 1.2 of 
\cite{ber-mar:2006} since when $\cala(w,N)$ is Koszul one has $\cala^!_{\nu_N(n)}=\Ext^n_\cala(\mathbb K, \mathbb K)$ and the product of $\cala'(w,N)$ is essentially the Yoneda product (\cite{ber-mar:2006}, Proposition 3.1). Let us recall this theorem 1.2 of \cite{ber-mar:2006} which is an important result.

\begin{theorem}\label{BerM}
Let $\cala$ be a $N$-homogeneous algebra which is Koszul of finite global dimension. Then $\cala$ is Gorenstein if and only if the Yoneda algebra $E(\cala)=\Ext_\cala(\mathbb K,\mathbb K)$ is Frobenius.
\end{theorem}

As pointed out before this follows from Theorem \ref{KGD} and Theorem \ref{APR} by using the fact that one has $\cala'=E(\cala)$ whenever $\cala$ is Koszul.\\

\noindent \underbar{Remarks}.\\
 1) One sees that, with $D,N$ and $w$ as in Theorem \ref{APR}, one has two natural Frobenius algebras associated with $\cala(w,N)$. The first one is the algebra $\calf(w,N)=\cala^!/\cali$ of Theorem \ref{MOD} the other one is the algebra $\calf'(w,N)=\cala'/\cali'$ where 
\[
\cali'=\{ y\in \cala' \vert \omega_w(x\bullet y)=0,\ \ \forall x\in \cala'\}
\]
is a two-sided ideal since $\sigma_w$ induces an automorphism of $\cala'$ satisfying $\omega_w(x\bullet y)=\omega_w(\sigma_w(y)\bullet x)$. These two Frobenius algebras coincide for $N=2$ but are different for $N>2$.\\ 
2) $D,N$ and $w$ being as in Theorem \ref{APR}, it is tempting in view of Proposition \ref{Eq3R} to say that $w$ is $D$-{\sl regular} whenever the equivalent conditions (a) and (b) are satisfied. In fact Condition (a) contains several nondegeneracy conditions. This notion involves both $D$ and $N$ as above.

\subsection{Examples}

Of course one has already all the examples of Section 3. Let us give two quadratic examples and a class of $N$-homogeneous examples.\\

\noindent (a) {\sl The extended 4-dimensional Sklyanin algebra} \cite{ac-mdv:2002a}, \cite{ac-mdv:2003}, \cite{ac-mdv:2008}.
In connection with a problem of $K$-homology, the following quadratic algebra $\cala_\ug$ has been been introduced in \cite{ac-mdv:2002a} and analyzed in details in \cite{ac-mdv:2003}, \cite{ac-mdv:2008}. The algebra $\cala_\ug$ is the quadratic algebra generated by 4 elements $x^\lambda$ ($\lambda\in \{0,1,2,3\}$) with relations
\begin{equation}
\cos (\varphi_0-\varphi_k)[x^0,x^k]=i\sin (\varphi_\ell -\varphi_m)[x^\ell,x^m]_+
\label{ASk1}
\end{equation}
\begin{equation}
\cos(\varphi_\ell-\varphi_m)[x^\ell,x^m]=i\sin(\varphi_0-\varphi_k)[x^0,x^k]_+
\label{ASk2}
\end{equation}
for any cyclic permutation $(k,\ell,m)$ of (1,2,3). The parameter $\ug$ is the element
$\ug=\left(e^{i(\varphi_1-\varphi_0)},e^{i(\varphi_2-\varphi_0)},e^{i(\varphi_3-\varphi_0)}\right)$ of $T^3$. Thus there are a priori 3 scalar parameters $\varphi_1-\varphi_0,\varphi_2-\varphi_0$ and $\varphi_3-\varphi_0$. However for generic values of these parameters one can show that $\cala_\ug$ only depends on two scalar parameters and that then by an appropriate linear change of generators it reduces to the 4-dimensional Sklyanin algebra introduced in \cite{skl:1982} and studied in \cite{smi-sta:1992} from the point of view of general regularity.\\
The algebra $\cala_\ug$ is Koszul of global dimension $D=4$ and is Gorenstein whenever none of the 6 relations (\ref{ASk1}), (\ref{ASk2}) becomes trivial and one then has the nontrivial Hochschild cycle (in $Z(\cala,\cala)$)
\[
\begin{array}{lll}
w=\tilde{ch}_{\frac{3}{2}}(U_\ug)& = &-\sum_{\alpha, \beta,\gamma,\delta}\varepsilon_{\alpha \beta \gamma\delta} \cos (\varphi_\alpha-\varphi_\beta+\varphi_\gamma-\varphi_\delta)x^\alpha \otimes x^\beta \otimes x^\gamma \otimes x^\delta \\
\\
& + & i\sum_{\mu,\nu}\sin(2(\varphi_\mu-\varphi_\nu))x^\mu\otimes x^\nu \otimes x^\mu \otimes x^\nu
\end{array}
\]
which defines a 4-linear form on $\mathbb K^4$ which is preregular with 
\[
Q_w=-\bbbone
\]
i.e. $w$ is graded-cyclic. One verifies that one has then $\cala_\ug=\cala(w,2)$ and that $\bbbone\otimes w$ is a Hochschild 4-cycle, i.e. $\bbbone \otimes w \in Z_4(\cala,\cala)$.\\

\noindent (b) {\sl The $q$-deformed $D$-dimensional polynomial algebra}. This is the algebra $\cala$ generated by $D$ elements $x^\lambda$ ($\lambda\in \{1,\dots,D\}$) with relations
\begin{equation}
x^\mu x^\nu=q^{\mu\nu} x^\nu x^\mu
\label{qdefD}
\end{equation}
for $\mu,\nu\in \{1,\dots,D\}$ where the $q^{\mu\nu}\in \mathbb K$ satisfy
\begin{equation}
q^{\mu\nu}q^{\nu\mu}=1,\ \ \ q^{\lambda\lambda}=1
\label{qrel}
\end{equation}
for any $\lambda,\mu,\nu\in \{1,\dots,D\}$.\\
This algebra is Koszul of global dimension $D$ and Gorenstein. One has $\cala=\cala(w,2)$ with
\begin{equation}
w=\sum_{\pi\in \fracS_D} \chi(\pi)x^{\pi(1)}\otimes \dots \otimes x^{\pi(D)}
\label{wqdefD}
\end{equation}
where $\fracS_D$ is the group of permutations of $\{1,\dots,D\}$ and where $\chi:\fracS_D\rightarrow \mathbb K$ is given by $\chi(\pi)=\prod_{(\mu\nu)}(-q^{\mu\nu})$ with
$\Pi_{(\mu\nu)}$ corresponding to the standard embedding
\[
\fracS_D\hookrightarrow\{\prod_{(\mu\nu)} b^{\mu\nu}, \mu<\nu\} \subset \fracB_D
\]
of $\fracS_D$ into the group of braids $\fracB_D$.\\
One has then
\begin{equation}
(Q_w)^\mu_\nu = \left(\prod_{\lambda\not= \mu}(-q^{\lambda\mu})\right)\delta^\mu_\nu
\label{QqdefD}
\end{equation}
for the matrix element of the corresponding $Q_w\in GL(D,\mathbb K)$.\\

\noindent (c) {\sl Precommutative examples} \cite{ber:2001a},\cite{ber-mar:2006}. Let the integers $g$ and $N$ be such that $g\geq N\geq 2$ and let $\varepsilon$ be the completely antisymmetric $g$-linear form on $\mathbb K^g$ with $\varepsilon(e_1,\dots,e_g)=1$. Consider the $N$-homogeneous algebra $\cala=\cala(\varepsilon, N)$ i.e. the algebra generated by $g$ elements $x^\lambda$ ($\lambda\in \{1,\dots,g\}$) with the relations
\[
\varepsilon_{\lambda_1\dots \lambda_{g-N}\ \mu_1\dots \mu_N}x^{\mu_1}\dots x^{\mu_N}=0
\]
where $\varepsilon_{\lambda_1\dots \lambda_g}=\varepsilon(e_{\lambda_1},\dots,e_{\lambda_g})$. It is clear that $\varepsilon$ is preregular with
\[
Q_\varepsilon= (-1)^{g-1}\bbbone
\]
as associated element of $GL(g,\mathbb K)$.\\
It was shown in \cite{ber:2001a} where this algebra was introduced that $\cala(\varepsilon,N)$ is a Koszul algebra of finite global dimension and it was shown in 
\cite{ber-mar:2006} that it is Gorenstein if and only if either $N=2$ or $N>2$ and $g=Np+1$ for some integer $p\geq 1$. For $N=2$ this reduces to the algebra polynomial functions on $\mathbb K^g$ while for $N>2$ and $g=Np+1$ this is a regular algebra of global dimension $D=2p+1$. In the latter case, the ideal $\cali$ of Theorem \ref{MOD} is generated by the quadratic elements $\alpha\beta + \beta \alpha$ of $\cala(\varepsilon,N)^!$
so that the quotient Frobenius algebra $\calf (\varepsilon, N)=\cala^!/\cali$ reduces to the exterior algebra $\wedge\mathbb K^g$ which is precisely the Koszul dual algebra of the quadratic algebra of polynomial functions on $\mathbb K^g$. Thus by this process one recovers the quadratic relations implying the original $N$-homogeneous ones.\\
Notice that for $N>2$ the algebra $\cala(\varepsilon, N)$ has exponential growth \cite{ber:2001a}.\\
In \cite{hai-kri-lor:2008} a twisted version of this example associated with a Hecke symmetry is introduced and analyzed with similar results. This paper \cite{hai-kri-lor:2008} contains even a super version of these examples. See also \cite{gur:1990} (and \cite{wam:1993}) for the quadratic case associated with a Hecke symmetry.

\noindent \underbar{Remark}. In contrast to the previous example for $N>2$, in the cases of the Yang-Mills algebra and the super Yang-Mills algebra the ideal $\cali$ of Theorem \ref{MOD} vanishes, that is the Koszul duals are then Frobenius. The reason is that in these cases the 3-regular multilinear forms (4-linear) $w$ given respectively by (\ref{wYM}) and (\ref{wSYM}) satisfy the stronger condition (iii') of \S 3.1.

\subsection{Classical limit versus infinitesimal preregularity}

We now consider perturbations of the algebra $\mathbb K[x^1,\dots,x^g]$ of polynomial functions on $\mathbb K^g$. More precisely one has $\mathbb K[x^1,\dots,x^g]=\cala(\varepsilon,2)$ where $\varepsilon$ is the $g$-linear form on $\mathbb K^g$ which is completely antisymmetric with $\varepsilon_{1\ 2\dots g}=1$, where $\varepsilon_{\lambda_1\dots\lambda_g}=\varepsilon(e_{\lambda_1},\dots, e_{\lambda_g})$ are the components of $\varepsilon$ in the canonical basis $(e_\lambda)$ of $\mathbb K^g$. Let $w_t$ be a 1-parameter family of preregular $g$-linear forms on $\mathbb K^g$ with $w_0=\varepsilon$ and let us investigate what happens formally at first order in $t$. One writes
\begin{equation}
\left\{
\begin{array}{l}
w_t=\varepsilon + t \dot w + o(t^2)\\
Q_{w_t}=(-1)^{g-1}\bbbone + t \dot Q + o(t^2)
\end{array}
\right.
\label{1-prer}
\end{equation}
and the first order $Q_{w_t}$-cyclicity reads
\begin{equation}
\dot W_{\lambda_1\dots\lambda_g}=\dot Q^\lambda_{\lambda_g} \varepsilon_{\lambda\lambda_1\dots \lambda_{g-1}}+(-1)^{g-1} \dot W_{\lambda_g\lambda_1\dots \lambda_{g-1}}
\label{1cycl}
\end{equation}
with $\dot W_{\lambda_1\dots \lambda_g}=\dot w(e_{\lambda_1},\dots, e_{\lambda_g})$.
This equation implies
\begin{equation}
\tr (\dot Q)=\dot Q^\lambda_\lambda=0
\label{detQ}
\end{equation}
which suggests $\det (Q_{w_t})=1$ for a finite version. So a natural question is the following : Does a quadratic AS-regular algebra $\cala(w,2)$ is such that $\det(Q_w)=1$? By looking at Example (c) of \S 3.3, one can see that the answer is no. Notice however that the quadratic AS-algebra of type $E$ is isolated.


\end{document}